\documentclass[final,1p,times]{elsarticle}
\usepackage{epsfig,amsmath,amssymb,amsthm,color,amsfonts,epstopdf,mathrsfs}
\biboptions{numbers,sort&compress}

\everymath{\displaystyle}
\def\l{\left}
\def\r{\right}

\def\RR{{\mathbb{R}}}

\def\DD{{\cal{D}}}
\def\OO{{\cal{O}}}
\def\SS{{\cal{S}}}
\def\R#1{$(\ref{#1})$}

\newcommand{\bb}[1]{\begin{equation}\label{#1}}
\newcommand{\ee}{\end{equation}}
\newcommand{\bbb}{\begin{eqnarray}}
\newcommand{\eee}{\end{eqnarray}}
\newcommand{\bbbb}{\begin{eqnarray*}}
\newcommand{\eeee}{\end{eqnarray*}}
\newcommand{\tT}{\intercal}
\newcommand{\nnn}{\nonumber}

\definecolor{green1}{rgb}{0.1,0.5,0.0}

\newtheorem{thm}{Theorem}
\newtheorem{lemma}{Lemma}
\newtheorem{cor}{Corollary}
\theoremstyle{remark}
\newtheorem{rem}{Remark}
\theoremstyle{define}
\newtheorem{define}{Definition}

\newcommand{\clearallnum}{
    \numberwithin{equation}{section} \setcounter{equation}{0}
    \numberwithin{thm}{section} \setcounter{thm}{0}
    \numberwithin{lemma}{section} \setcounter{lemma}{0}
    \numberwithin{cor}{section} \setcounter{cor}{0}
    \numberwithin{rem}{section} \setcounter{rem}{0}
    \numberwithin{define}{section} \setcounter{define}{0}}

\journal{Applied Mathematics and Computation}

\begin{document}

\begin{frontmatter}



\title{Numerical solution of degenerate stochastic Kawarada equations via a semi-discretized approach}


\author{Joshua L. Padgett\corref{cor}\fnref{label}}
\ead{josh\_ padgett@baylor.edu}
\fntext[label]{Corresponding Author}
\author{Qin Sheng}
\ead{qin\_ sheng@baylor.edu}

\address{Department of Mathematics and
Center for Astrophysics, Space Physics and Engineering Research,
Baylor University, Waco, TX 76798-7328, USA}

\begin{abstract}
The numerical solution of a highly nonlinear two-dimensional degenerate stochastic Kawarada equation 
is investigated. A semi-discretized approximation in space is comprised on arbitrary nonuniform grids. 
Exponential splitting strategies are then applied to advance solutions of the semi-discretized scheme over adaptive 
grids in time. It is shown that key quenching solution features including the positivity and monotonicity are well
preserved under modest restrictions. The numerical stability of the underlying splitting method is also maintained 
without any additional restriction. Computational experiments are provided to not only illustrate our results, but also 
provide further insights into the global nonlinear convergence of the numerical solution.

\end{abstract}


\end{frontmatter}


\section{Introduction} \clearallnum
Kawarada partial differential equations arise in the modeling of highly sophisticated, yet important, natural phenomena
where singularities may develop as the solution evolves in time. Such singularities often result from the energy of a
system concentrating and approaching its activation criterion \cite{Poin,Bebernes_89}. 

Consider a typical solid-fuel ignition process. If the combustion chamber with fuel and air are appropriately mixed, the 
temperature in the chamber may increase monotonically until a certain critical value is reached. However, rates of
such temperature changes can occur in a nonlinear manner throughout the media. This nonuniform distribution of heat 
may result in high temperatures being extremely localized within the combustion chamber and lead to an
ignition once the peak temperature reaches a certain threshold \cite{Poin}. This phenomenon is carefully 
characterized by a Kawarada model in which temporal derivatives of the solution may grow at an explosive rate, 
while the solution itself remains bounded \cite{Chan2,Kawa,Josh1,Josh2}. This strong nonlinear singularity, which is 
referred to as a {\em quenching singularity}, has been the backbone of the modeling equations. Applications of 
Kawarada equations can be found in fuel combustion simulations, enhanced thermionic emission optimization processes, 
electric current transients in polarized ionic chambers, and enzyme kinematics.

Let $\DD=(-a,a)\times (-b,b)\subset\RR^2$ be an idealized combustion chamber, $\partial\DD=\bar{\DD}\setminus\DD$ be
its boundary, and $\Omega=\DD\times (t_0,T_q),~ \SS=\partial\DD\times(t_0,T_q),$ where $a,b>0,~0\le t_0<T_q<\infty$ 
with $T_q$ being the terminal time. Further, let $c>0$ be the ignition temperature. If $u(x,y,t)$ denotes the temperature 
distribution within the chamber, we have the following 2D model generalized from the original Kawarada configuration \cite{Kawa}:
\bbb
&& \mathscr{L}u=f(\varepsilon,u),~~~(x,y,t)\in\Omega,\label{b1}\\
&& u(x,y,t)=0,~~~(x,y,t)\in\SS,\label{b2}\\
&& u(x,y,t_0)=u_0(x,y)\ll c,~~~(x,y)\in \DD,\label{b3}
\eee
where
$$\mathscr{L} \mathrel{\mathop:}= \sigma(x,y)\partial_t-\Delta,$$
and $\Delta$ is the standard 2D Laplacian. The degeneracy function $\sigma(x,y)\geq 0$ for $x\in\bar{\DD}$ with equality 
occurring only for $(x,y)\in\partial\DD_0\subseteq\partial\DD.$
The nonlinear reaction term, $f(\varepsilon,u)$ with a bounded stochastic influence $\varepsilon(x,y),$  
satisfies the following properties
$$f(\varepsilon,0)=f_0(\varepsilon)>0,~~f_u(\varepsilon,u)>0,~~\lim_{u\to c^-}f(\varepsilon,u)=\infty,~~
\int_0^c f(\varepsilon,u)\,du = \infty$$
for $u\in[0,c).$ The stochastic influence is characterized by
$$0<f^\varepsilon_{\min} \le \min_{\varepsilon}\{f(\varepsilon,u)\}\quad\mbox{and}\quad \max_{\varepsilon}\{f(\varepsilon,u)\} \le f^{\varepsilon}_{\max} \ll \infty$$
for each $u\in[0,c).$
The solution
$u$ of \R{b1}-\R{b3} is said to {\em quench\/} if there exists $T_q<\infty$ such that
\bb{a1}
\sup\l\{u_t(x,y,t)\,:\,(x,y)\in {\DD}\r\}\rightarrow\infty~\mbox{as}~t\rightarrow
T_q^{-}.
\ee
$T_q$ is then referred to as the {\em quenching time\/} \cite{Chan2,Acker2,Acker1}.
It has been shown that a necessary condition for quenching to occur with the above conditions placed on the reaction term is
\bb{a2}
\max\l\{u(x,y,t)\,:\,(x,y)\in \bar{\DD}\r\}\rightarrow c^{-}~\mbox{as}~t\rightarrow T_q^{-}.
\ee
We may note that quenching is a mathematical interpretation of the fuel ignition/combustion process
in physics. Since $T_q<\infty$ only when certain spatial references, such as the size and shape of $\DD,$ 
reach their critical limits, a domain $\Omega_*$ is called a {\em critical domain} of \R{b1}-\R{b3} if 
the solution of \R{b1}-\R{b3} exists globally as $\Omega\subseteq\Omega_*$ and \R{a2} occurs 
as $\Omega_*\subseteq\Omega.$ We note that such critical domains are not unique in 
multi-dimensional circumstances, since both the {\em size} and {\em shape} of $\Omega$ affect 
combustion \cite{Chan3,Beau1,Beau2,Sheng5}. 

Highly effective and efficient algorithms for solving Kawarada equations such as \R{b1}-\R{b3} have been 
difficult to obtain due to the highly nonlinear mechanism. There has been much effort to develop 
algorithms that may accurately predict the ignition location and time, solution profiles,
and critical domain characteristics, while reasonably conserving the solution positivity and monotonicity \cite{Sheng5,Sheng3,Sheng4,Sheng6,Chan1,Chan4}. While explorations have 
been carried out on both uniform and nonuniform grids, as well as via adaptations, rigorous
analysis of different schemes involving nonuniform grids has been incomplete. 
To improve computational efficiency, splitting strategies are introduced and incorporated with adaptive
mesh structures
\cite{Sheng3,Lang2,Beau2}. The success of such combinations has also become a key motivation for
the study to be presented in this paper.

Rigorous mathematical analysis is implemented for a highly vibrant semi-discretization
oriented method for solving \R{b1}-\R{b3} in our study.
In the next section, we verify carefully the order of accuracy of the scheme proposed. In Sections 3 and 4, we prove 
that the solution of the semi-discretized system, as well as the corresponding scheme, preserve the 
expected solution positivity and monotonicity under reasonable restrictions. Section 5 is devoted to standard and 
reinforced stability 
analysis of the numerical method. Experimental results to validate our analysis and explore global nonlinear convergence 
of the stochastic singular Kawarada solutions are presented in Section 6. Finally, Section 7 is tailored for concluding our
remarks, concerns, and forthcoming endeavors.

\section{Derivation of Scheme} \clearallnum

Without loss of generality, we set $c=1.$ Utilizing transformations $\tilde{x}=x/a,~\tilde{y}=y/b,$ and reusing the 
original variables for notation simplicity, we may reformulate \R{b1}-\R{b3} as
\bbb
&&u_t=\frac{1}{a^2\sigma}u_{xx} + \frac{1}{b^2\sigma}u_{yy} +g(\varepsilon,u),~~(x,y,t)\in\Omega,~~~~~~\label{c1}\\
&&u(x,y,t)=0,~~(x,y,t)\in \SS,\label{c2}\\
&&u(x,y,t_0)=u_0(x,y)\ll 1,~~(x,y)\in \DD,\label{c2b}
\eee
where $\DD=(-1,1)\times (-1,1)\subset\RR^2$ and $g(\varepsilon,u) = f(\varepsilon,u)/\sigma(x,y).$
For fixed $N_x, N_y\gg 1$ we define
\bbbb
\DD_h&=&\l\{(x_i,y_j)\,:\, 0\leq i\leq N_x+1; ~0\leq j\leq N_y+1\r\},\\
\DD_h^{\circ}&=&\l\{(x_i,y_j)\,:\, 0< i< N_x+1; ~0< j< N_y+1\r\},
\eeee
where $x_0=y_0=-1,~x_{i+1}-x_i=h_{x,i}>0,~y_{j+1}-y_j=h_{y,j} >0$ for $0\leq i\leq N_x,~0\leq j\leq N_y,$ 
and $x_{N_x+1}=y_{N_y+1}=1.$ Let $u_{i,j}(t)$ denote an approximation of the solution of \R{c1}-\R{c2b} at  
$(x_i,y_j,t)\in\DD_h\times [t_0, T_q).$ We consider the following nonuniform finite difference formulas \cite{Sheng3},
\bbb
\l.\frac{\partial^2u}{\partial x^2}\r|_{(x_i,y_j,t)} &\approx& \frac{2u_{i-1,j}(t)}{h_{x,i-1}(h_{x,i-1}+h_{x,i})}-\frac{2u_{i,j}(t)}{h_{x,i-1}h_{x,i}}+\frac{2u_{i+1,j}(t)}{h_{x,i}(h_{x,i-1}+h_{x,i})},~~~~~~\label{xder}\\
\l.\frac{\partial^2u}{\partial y^2}\r|_{(x_i,y_j,t)} &\approx& \frac{2u_{i,j-1}(t)}{h_{y,j-1}(h_{y,j-1}+h_{y,j})}-\frac{2u_{i,j}(t)}{h_{y,j-1}h_{y,j}}+\frac{2u_{i,j+1}(t)}{h_{y,j}(h_{y,j-1}+h_{y,j})},~~~~~~\label{yder}
\eee
for $(x_i,y_j)\in\DD_h^{\circ}.$ Further, for a standard lexicographical ordering, we denote
$v(t)=(u_{1,1},u_{2,1},\dots,$ $u_{N_x,1},u_{1,2},\dots,u_{N_x,N_y})^{\tT}\in\RR^{N_xN_y},$
and let $g(\varepsilon,v)$ be a suitable discretization of the nonlinear stochastic reaction term. 
We derive from from \R{c1}-\R{c2b} our semi-discretized system
\bbb
v'(t)&=& Mv(t)+g(\varepsilon,v(t)),~~~t_0<t<T_q,\label{c8}\\
v(t_0)&=&v_0.\label{c8b}
\eee
Denote $h_{\min} = \textstyle\min_{1\leq i\leq N_x;1\leq j\leq N_y}\{h_{x,i-1}h_{x,i},h_{y,j-1}h_{y,j}\}$ and
$$M_x=\frac{1}{a^2}B(I_{N_y}\otimes T_x),~M_y=\frac{1}{b^2}B(T_y
\otimes I_{N_x}),~M=M_x+M_y,$$
where $\otimes$ stands for the Kronecker product,
$I_{N_{\alpha}}\in\RR^{N_{\alpha}\times N_{\alpha}},~\alpha=x,y,$ are identity matrices. Further, let
\bbbb
B&=&\mbox{diag}\l(\sigma_{1,1}^{-1},\sigma_{2,1}^{-1},\dots,\sigma_{N_x,1}^{-1},\sigma_{1,2}^{-1},\dots,\sigma_{N_x,N_y}^{-1}\r)\in\RR^{N_xN_y\times N_xN_y},\\
T_{\alpha}&=&\mbox{tridiag}\l(h_{\min}^{-1}l_{\alpha,i},h_{\min}^{-1}m_{\alpha,i},h_{\min}^{-1}n_{\alpha,i}\r)\in
\RR^{N_{\alpha}\times N_{\alpha}},~~~\alpha=x,y,
\eeee
and for the above,
\bbbb
l_{\alpha,j} &=& \frac{2h_{\min}}{h_{\alpha,j}(h_{\alpha,j}+h_{\alpha,j+1})},~~~n_{\alpha,j} ~=~ \frac{2h_{\min}}{h_{\alpha,j}(h_{\alpha,j-1}
+h_{\alpha,j})},~~~j = 1,\dots,N_{\alpha}-1,\\
m_{\alpha,j} &=& -\frac{2h_{\min}}{h_{\alpha,j-1}h_{\alpha,j}},~~~j=1,\dots,N_{\alpha};~\alpha=x,y.
\eeee
The formal solution of \R{c8}, \R{c8b} is
\bb{rs1}
v(t)=E(tM)v_0+\int_{t_0}^tE((t-s)M)g(\varepsilon,v(s))\,ds,~~~t_0<t<T_q,
\ee
where $E(\cdot)$ is the corresponding matrix exponential \cite{Sheng3,Golub,Iserles}. 
An aim of this paper is to investigate a highly efficient numerical procedure for the problem \R{c8}-\R{rs1}. 

To this end, we consider the following realization of \R{rs1} on variable temporal grids:
\bb{xxx}
v_{k+1}=E(\tau_k M)v_k+\int_{t_k}^{t_{k+1}}E((t_{k+1}-s)M)g(\varepsilon,v(s))\,ds,~~~k=0,1,\dots
\ee
where $v_k$ and $v_{k+1}$ are approximations to $v(t_k)$ and $v(t_{k+1}),$ respectively, $v_0$ is 
the initial vector, $t_k = t_0+\textstyle\sum_{j=0}^{k-1}\tau_j,$ $0<\tau_k\ll 1,$ $k=0,1,2,\dots,$ and 
$\{\tau_k\}_{k\ge 0}$ is a set of adaptive temporal steps. %% \cite{Fur,Sheng5,Beau3,}. 
In fact, a rapid temporal step refinement mechanism must be triggered robotically as $t_k\to T_q^-,$ since $u_t$ may grow exponentially and affect computational accuracy near quenching points
\cite{Poin,Chan2,Levine}. On the other hand, an overly refined mesh should also be avoided 
\cite{Fur,Sheng8}. In light of these concerns, a standard arc-length monitoring function may be implemented to handle 
the strong nonlinearity at quenching \cite{Sheng3,Sheng4,Lang2,Fur,Sheng5,Beau3}.

Further, to avoid the need for nonlinear iterations, in computations we replace $g(\varepsilon,v_{k+1})$ by
$g(\varepsilon,w_{k}),$ where $w_{k}$ is a suitable approximation to $v_{k+1}$ with a proper order of accuracy \cite{Beau1,Josh1,Josh2}.
However, no such approximation is utilized in the following analysis.

We have the following weighted approximation.
\begin{lemma}
Assume that $f(\varepsilon,u)\in C_c^2[0,1)$ and $\tau_k\rightarrow 0^+.$ Then for $\theta\in[0,1]$ we have
$$\int_{t_k}^{t_{k+1}}E((t_{k+1}-s)M)g(\varepsilon,v(s))\,ds - M^{-1}[E(\tau_kM)-I][\theta g(\varepsilon,v_k)+(1-\theta)g(\varepsilon,v_{k+1})]=\OO\l(\tau_k^p\r),$$ 
where $p=2$ when $\theta = 1/2,$ and $p=1$ otherwise.  
\end{lemma}

\begin{proof}
We begin with the following difference
$$\eta_k=\int_{t_k}^{t_{k+1}}E((t_{k+1}-s)M)g(\varepsilon,v(s))\,ds - \int_{t_k}^{t_{k+1}}E((t_{k+1}-s)M)[\theta g(\varepsilon,v_k)+(1-\theta)g(\varepsilon,v_{k+1})]\,ds.$$
Replace the first integral by a trapezoidal rule and integrate the second directly to obtain
$$\eta_k=\frac{\tau_k}{2}\l[g(\varepsilon,v_{k+1})+E(\tau_kM)g(\varepsilon,v_k)\r]
- M^{-1}[E(\tau_kM)-I][\theta g(\varepsilon,v_k)+(1-\theta)g(\varepsilon,v_{k+1})]+\OO\l(\tau_k^3\r),$$
as $\tau_k\rightarrow 0^+.$ 
%We now consider local expansions of each of the terms
%$$g(\varepsilon,v_{k+1}) = g(\varepsilon,v_k)+\tau_kg'(\varepsilon,v_k)+\frac{\tau_k^2}{2}%g''(\varepsilon,v_k)+\OO\l(\tau_k^3\r)$$
%$$E(\tau_kM)g(\varepsilon,v_k) = g(\varepsilon,v_k)+\tau_kMg(\varepsilon,v_k)+\frac{\tau_k^2}{2}M^2g(\varepsilon,v_k)+\OO\l(\tau_k^3\r)$$
%$$M^{-1}[E(\tau_kM)-I][\theta g(\varepsilon,v_k)] = \tau_k\theta g(\varepsilon,v_k)+\frac{\tau_k^2\theta}{2}Mg(\varepsilon,v_k)+\OO\l(\tau_k^3\r)$$
%$$M^{-1}[E(\tau_kM)-I][(1-\theta) g(\varepsilon,v_{k+1})] = (1-\theta)\l[\tau_k+\frac{\tau_k^2}{2}M\r]g(\varepsilon,v_k)+\tau_k^2(1-\theta)g'(\varepsilon,v_k)+\OO\l(\tau_k^3\r)$$
Consider local expansions of each term. A direct comparison then leads to. %%yields
% \bbbb
% &&\frac{\tau_k}{2}\l[g(\varepsilon,v_{k+1})+E(\tau_kM)g(\varepsilon,v_k)\r]+\OO\l(\tau_k^3\r)\\
% &&~~~~~~~~~~~~~~~~~~~~~~~~~~~~~~~~~~~~~~~~ - M^{-1}[E(\tau_kM)-I][\theta g(\varepsilon,v_k)+(1-\theta)g(\varepsilon,v_{k+1})]\\
$$\eta_k=\frac{\tau_k^2}{2}g_{v_k}(\varepsilon,v_k) - (1-\theta)\tau_k^2g_{v_k}(\varepsilon,v_k)+
\OO\l(\tau_k^3\r),~~~\tau_k\rightarrow 0^+,$$
% \eeee
and apparently, the approximation is first-order unless $\theta=1/2,$ in which case the approximation is
second-order.
\end{proof}

Let $\Phi_{x,k} \equiv E(\tau_kM_x/2)$ and $\Phi_{y,k} \equiv E(\tau_kM_y).$
We further denote $\Phi_k \equiv \Phi_{x,k}\Phi_{y,k}\Phi_{x,k}$ and
$F(\varepsilon,v_k,v_{k+1},\theta)\equiv \theta g(\varepsilon,v_k)+(1-\theta)g(\varepsilon,v_{k+1}).$ 
Thus, from \R{xxx} and Lemma 2.1 we acquire that
%\bbb
%v_{k+1}&\approx &E(\tau_k M)v_k+\int_{t_k}^{t_{k+1}}E((t_{k+1}-s)M)g(\varepsilon,v_k)ds\nnn\\
%&=&E(\tau_kM)v_k+M^{-1}\l[E(\tau_kM)-I\r]g(\varepsilon,v_k)%\label{A1}
%\eee
%\bbb
%v_{k+1}&\approx &E(\tau_k M)v_k+\int_{t_k}^{t_{k+1}}E((t_{k+1}-s)M)g(\varepsilon,v_{k+1})ds\nnn\\
%&=&E(\tau_kM)v_k+M^{-1}\l[E(\tau_kM)-I\r]g(\varepsilon,v_{k+1})%\label{A2}
%\eee
\bbb
v_{k+1}&=&E(\tau_k M)v_k+\int_{t_k}^{t_{k+1}}E((t_{k+1}-s)M)F(\varepsilon,v_k,v_{k+1},\theta)\,ds +\OO\l(\tau_k^{p+1}\r)\nnn\\
%&=&E(\tau_kM)v_k+M^{-1}\l[E(\tau_kM)-I\r]F(\varepsilon,v_k,v_{k+1},\theta)+\OO\l(\tau_k^{p+1}\r)\nnn\\
&= &\Phi_k v_k+ M^{-1}\l[\Phi_k-I\r]F(\varepsilon,v_k,v_{k+1},\theta)+\OO\l(\tau_k^{p+1}\r),~~~
\tau_k\rightarrow 0^+,\label{A3}
\eee
where $\theta\in[0,1]$ is a fixed parameter, and $p\in\{1,2\}$ depending upon the choice of $\theta.$ We note that the stochastic influences are completely contained in the vector-valued function $F(\varepsilon,v_k,v_{k+1},\theta).$

Before proceeding, we would like to note that, in general, \R{xder}, \R{yder} are first order approximations \cite{Sheng3,Iserles}. The following lemma demonstrates that in certain cases the accuracy can be increased to second-order.

\begin{lemma}
If the nonuniform mesh $\DD_h$ is the image of a uniform mesh on $\DD$ via a smooth mapping, then \R{xder} and \R{yder} are in fact second-order.
\end{lemma}

\begin{proof}
We consider the case for \R{xder} and note that the result for \R{yder} follows in a similar manner. Assume that $u\in C^4(\DD)$ and $\DD_h$ is generated via a smooth mapping $g^{(x)}$ defined as
\bb{1111}
x_i = g^{(x)}(\omega^{(x)}_i),\ i=0,\ldots,N+1
\ee
where $\omega^{(x)}_i = -1+ik_1,\ i=0,\ldots,N+1,$ with $k=2/(N+1).$ Thus, $u(x,y,t) = u(g^{(x)}(\omega^{(x)}),y,t).$ We now consider the local truncation error, $\psi^{(x)}_{i,j}(t),$ of \R{xder} at a point $(x_i,y_j,t)\in\DD_h\times[t_0,T_q):$
\bbb
\psi^{(x)}_{i,j}(t) &=& \frac{\partial^2 u}{\partial x^2}\bigg|_{(x_i,y_j)} -\l[\frac{2u_{i-1,j}(t)}{h_{x,i-1}(h_{x,i-1}+h_{x,i})}-\frac{2u_{i,j}(t)}{h_{x,i-1}h_{x,i}}+\frac{2u_{i+1,j}(t)}{h_{x,i}(h_{x,i-1}+h_{x,i})}\r]\nnn\\
& = & 2\frac{\partial^3 u}{\partial x^3}\bigg|_{(x_i,y_j)}(h_{x,i-1}-h_{x,i}) + \OO\l(\frac{h_{x,i}^3}{h_{x,i-1}+h_{x,i}}\r)\nnn\\
& = & 2\frac{\partial^3 u}{\partial x^3}\bigg|_{(x_i,y_j)}[(x_i-x_{i-1})-(x_{i+1}-x_i)] + \OO(\max\{h_{x,i-1}^2,h_{x,i}^2\}).\label{111a}
\eee
For general nonuniform meshes, \R{111a} is first-order, however, due to \R{1111}, we have
\bb{1122}
(x_i-x_{i-1})-(x_{i+1}-x_i) = -k_1^2\frac{x_{i-1}-2x_i+x_{i+1}}{k_1^2} = -k_1^2\frac{d^2g^{(x)}}{d(\omega^{(x)})^2}\bigg|_{\omega^{(x)}_i} + \OO(k_1^2).
\ee
Note that for $\xi_i = \epsilon \omega_i + (1-\epsilon)\omega_{i-1},\ \epsilon\in[0,1],$ we have
$$h_{x,i-1} = x_i-x_{i-1} = g^{(x)}(\omega_i) - g^{(x)}(\omega_{i-1}) = \frac{\partial g^{(x)}}{\partial\omega}\bigg|_{\xi_i}k_1$$
which gives that $k_1 = (g^{(x)}_\omega(\xi_i))^{-1}h_{x,i-1} = \OO(h_{x,i-1}).$ A similar process also gives $k_1 = \OO(h_{x,i}).$ These facts combined with \R{111a} and \R{1122} gives the desired result.
\end{proof}

\section{Positivity} \clearallnum

It has been shown that solutions to \R{b1}-\R{b3} remain positive within their domains of 
existence; {\em i.e}, for $t\in[t_0,T_q)$ \cite{Chan2,Acker2,Acker1,Levine}. This is also a key
characteristics of a healthy physical solution no matter whether an ignition, or quenching,
occurs \cite{Poin,Beau3}. The property must be properly preserved by \R{A3}, and likewise, 
by our semidiscretized problem \R{c8}, \R{c8b}. 
%To that end, let $\wedge$ be one of the operations $<,$ $\leq,$ $>,$ $\geq$ and $\alpha,\beta\in
%\RR^{K_1\times K_2}.$ We assume the following notations in subsequent discussions:

%\begin{enumerate}
%\item $\alpha\wedge\beta$ means $\alpha_{i,j}\wedge\beta_{i,j},~i=1,2,\ldots,K_1;\,j=1,2,\ldots,K_2;$

%\item $a\wedge\alpha$ means $a\wedge\alpha_{i,j},~~i=1,2,\ldots,K_1;\,j=1,2,\ldots,K_2,$
%for any given scalar $a.$
%\end{enumerate}
% We assume the following notations in continuing discussions:
% \begin{itemize}
% \item $\alpha\wedge\beta$ means $\alpha_{i,j}\wedge\beta_{i,j},~i=1,2,\ldots,K;\,j=1,2,\ldots,N;$
% \item $a\wedge\alpha$ means $a\wedge\alpha_{i,j},~~i=1,2,\ldots,K;\,j=1,2,\ldots,N,$
% for any given scalar $a.$
% \end{itemize}

All vector inequalities are to be interpreted componentwise in following investigations.

\begin{lemma}
Let $K>1$ and $A\in\RR^{K\times K}$ be nonsingular and nonnegative. Further, let
$\beta\in\RR^{K}$ be positive. Then $A \beta > 0.$
\end{lemma}

%\begin{proof}
%Let $A\in\RR^{K\times K}$ be nonsingular and nonnegative. In particular, this gives that each row of $A$
%must contain a nonzero element. Denote the rows of $A$ as $a_i$ for $i=1,\dots, K.$ Let $\beta\in\RR^K$
%be positive. Note that the product $A\beta = \l(\alpha_1,\dots,\alpha_K\r)^{\tT},$ where $\alpha_i = \l(a_i,\beta^{\tT}\r)$
%and $\l(\cdot,\cdot\r)$ is the standard Euclidean inner product. Since each $a_i$ has a positive element and
%$\beta$ is positive, it follows that $\alpha_i$ is positive for $i=1,\dots,K,$ giving the result.
%\end{proof}

\begin{proof}
The proof follows directly from the definitions.
\end{proof}

\begin{lemma}
The matrix $E(\tau_kM)$ is positive for all $\tau_k\ge 0.$
\end{lemma}

\begin{proof}
Note that $M = M_x+M_y$ is such that $M_{ij} \ge 0$ for $i\neq j.$ Let $N = M+\beta I,$ 
where $\beta = 4/h_{\min}\textstyle\min\{a^2,b^2\}.$ Then $N$ is a positive matrix and
$$E(\tau_kM) = E(\tau_k(N-\beta I)) = E(-\tau_k\beta I)E(\tau_kN).$$
Since $\beta I$ is a diagonal matrix, it follows that$E(-\tau_k\beta I)$ is a diagonal matrix, 
with $\l(E(-\tau_k\beta I)\r)_{ij} = e^{-\tau_k\beta} > 0$ for $i= j$ and $0$ elsewhere. 
Further, $E(\tau_k N)$ is positive by the definition of the matrix exponential and 
the fact that $N$ is positive. Thus, we have the result.
\end{proof}

%\begin{cor}
%The matrices $\Phi_{x,k}$ and $\Phi_{y,k}$ are positive for all $\tau_k\ge 0.$ Further, $\Phi_k$ is positive for all $\tau_k\ge 0.$
%\end{cor}

\begin{lemma}
The solution to \R{c8}, \R{c8b} is positive for all $t_0<t<T_q$ and $0\le u<1.$
\end{lemma} 

\begin{proof}
Note that the solution to \R{c8}, \R{c8b} is given by \R{rs1}, and that $E(tM)v_0$ is positive for all positive $v_0$ 
due to Lemmas 3.1 and 3.2. Similarly, $E((t-\tau)M)g(v(\tau))$ is positive since $g$ is a positive vector and by 
Lemma 3.1. Since the integral of positive vectors is positive, it follows naturally that $v(t)$ is positive.
\end{proof}

\begin{lemma}
The matrices $\Phi_{x,k}$ and $\Phi_{y,k}$ are positive for all $\tau_k\ge 0.$ Likewise, $\Phi_k$ is 
positive for all $\tau_k\ge 0.$
\end{lemma} 

\begin{proof}
The proof is similar to that of Lemma 3.2.
\end{proof}

We study the matrix $M^{-1}[\Phi_k-I]$ in the remainder of this section. To show that $M^{-1}[\Phi_k-I]$ is positive 
under appropriate restrictions of $\tau_k,$ we need the following definition and lemmas.

\begin{define}
Let $K>1$ and $A\in\mathbb{C}^{K\times K},$ $I_K\in\mathbb{C}^{K\times K}$ be the identity matrix. Further, 
let $\|\cdot\|$ be an induced matrix norm. Then the
associated logarithmic norm $\mu : \mathbb{C}^{K\times K}\to \RR$ of $A$ is defined as
$$\mu(A) = \lim_{h\to 0^+} \frac{\|I_K + hA\| - 1}{h}.$$
\end{define}

\begin{rem}
\cite{Golub} When the spectral norm is considered, we have $\mu_2(A) = \lambda_{\max}[(A+A^*)/2].$
\end{rem}

\begin{rem}
\cite{Desoer} When the maximal norm is considered, we have $\mu_\infty (A) = 
\textstyle\sup\{\mbox{Re}(a_{ii})+\textstyle\sum_{j\neq i}|a_{ij}|\}.$
\end{rem}

\begin{lemma}
%{\rm\cite{Golub}}  
For $t\ge 0$ we have $\|E(t A)\| \le E(t \mu(A))$ for any logarithmic norm, $\mu,$ 
and any induced matrix norm, $\|\cdot\|.$
\end{lemma}

\begin{proof}
Consider the linear initial value problem,
$$\begin{array}{ll}
x'(t)=Ax(t),& t>t_0,\\
x(t_0)=x_0>0,
\end{array}$$
which has the solution $x(t) = E(tA)x_0.$ Recall that the upper right Dini derivative 
for any function $x(t)$ with respect to $t$ is defined as
$$\l(D^+_t x\r)(t) = \limsup_{h\to 0^+}\frac{x(t+h)-x(t)}{h}.$$
Using this fact and taking the norm, we obtain the differential inequality
\bb{inq1}
D^+_t\|x(t)\| =  \limsup_{h\to 0^+}\frac{\|E(hA)x(t)\|-\|x(t)\|}{h} \leq \mu(A)\|x(t)\|,
\ee
where $\|\cdot\|$ is any induced matrix norm. Note that
\bb{inq2}
\|x(t)\| \le e^{t\mu(A)}\|x_0\|
\ee
is the solution to \R{inq1} \cite{Szarski}.
%which gives $$\frac{\|E(tA)x_0\|}{\|x_0\|} \le e^{t\mu(A)}.$$ 
Dividing both sides of \R{inq2} by $\|x_0\|\neq 0$ and taking the supremum over all $x_0$ 
such that $\|x_0\|\neq 0$ gives the desired result.
\end{proof}

\begin{lemma}
For any $\tau_k\ge 0,$ we have $\|\Phi_{x,k}\|_\infty \le 1$ and $\|\Phi_{y,k}\|_\infty \le 1.$
\end{lemma}

\begin{proof}
The proof follows immediately by applying Lemma 3.5 and then utilizing Remark 3.2 to compute 
the logarithmic norm of each matrix involved.
\end{proof}

Now, we wish to determine a bound for $\|M^{-1}\|_\infty$ under suitable conditions 
on our nonuniform spatial grids. To that end, we denote
\bbb
c_\alpha(h) &=& \min_{i=1,...,N_\alpha}\l\lbrace \frac{h_{\alpha,i}(2+h_{\alpha,i})}{h_{\alpha,i-1}(h_{\alpha,i-1}+h_{\alpha,i})} 
+ \frac{h_{\alpha,i+1}(h_{\alpha,i+1}-2)}{h_{\alpha,i}(h_{\alpha,i-1}+h_{\alpha,i})} + 2(h_{\alpha,i+1}-
h_{\alpha,i})\sum_{j=1}^ih_{\alpha,j}\r\rbrace,\nnn\\
k(h) &=& [a^{-2}c_x(h)+ b^{-2}c_y(h)]/4,\label{kval}
\eee
where $\alpha=x,y.$ Let $\sigma_{\max} = \textstyle\max_{1\leq i\leq N_x;1\leq j\leq N_y}\sigma_{i,j}.$
%$$c_y(h) = \min_{j=1,...,N_y}\l\lbrace \frac{h_{y,j}(2+h_{y,j})}{h_{y,j-1}(h_{y,j-1}+h_{y,j})} + 
%\frac{h_{y,j+1}(h_{y,j+1}-2)}{h_{y,j}(h_{y,j-1}+h_{y,j})} + 2(h_{y,j+1}-h_{y,j})\sum_{k=1}^jh_{y,k}\r\rbrace.$$ We also define 

\begin{lemma}
Assume that $k(h) \ge 1.$ Then $\|M^{-1}\|_\infty \le \sigma_{\max}/2.$
\end{lemma}

\begin{proof}
Let $\tilde{M} = B^{-1}M = a^{-2}(I_{N_y}\otimes T_x) + b^{-2}(T_y\otimes I_{N_x}).$ We have 
$$\|M^{-1}\|_\infty \le \|\tilde{M}^{-1}\|_\infty\|B^{-1}\|_\infty\le \sigma_{\max}\|\tilde{M}^{-1}\|_\infty.$$
Thus, we only need to consider the case with $\tilde{M}.$ To this end, we let $f_{i,j}$ be any mesh function 
defined on $\DD_h.$ Define $u_{i,j} = -\tilde{M}^{-1}f_{i,j},$ where $u_{i,j}$ solves the equation 
$-\tilde{M}u_{i,j} = f_{i,j}$ on $\DD_h$ and satisfies the homogeneous Dirichlet boundary conditions $u_{i,j}=0,~(x_i,y_j)\in\partial\DD_h.$ For $w(x,y) = [x^2 + y^2]/4,$ we define $w_{i,j}$ on $\overline{\DD}_h$ by
$$w_{i,j} = \frac{1}{4}\l[\l(-1+\sum_{k=1}^i h_{x,k}\r)^2+\l(-1+\sum_{k=1}^j h_{y,k}\r)^2\r].$$
Clearly, $\Delta w = w_{xx}+w_{yy} = 1.$ We now consider $\tilde{M}w_{i,j}$ by considering $\overline{\Delta}_{x,h}w_{i,j}\equiv 4(I_{N_y}\otimes T_x)w_{i,j}$ and $\overline{\Delta}_{y,h}w_{i,j}\equiv 4(T_y\otimes I_{N_x})w_{i,j}$ separately. Thus, we have
\bbbb
\overline{\Delta}_{x,h}w_{i,j} & = & l_{x,i-1}\l[\l(-1+\sum_{k=1}^{i-1} h_{x,k}\r)^2\r] + m_{x,i}\l[\l(-1+\sum_{k=1}^{i} h_{x,k}\r)^2\r]\\
%& = & l_{x,i-1}\l[\l(-1-h_{x,i}+\sum_{k=1}^{i} h_{x,k}\r)^2\r] + m_{x,i}\l[\l(-1+\sum_{k=1}^{i} h_{x,k}\r)^2\r]\\
&&+ n_{x,i}\l[\l(-1+\sum_{k=1}^{i+1} h_{x,k}\r)^2\r] + n_{x,i}\l[\l(-1+h_{x,i+1}+\sum_{k=1}^{i} h_{x,k}\r)^2\r]\\
%& = & l_{x,i-1}h_{x,i}\l(-2\sum_{k=1}^i h_{x,k} + 2 + h_{x,i}\r) + n_{x,i}h_{x,i+1}\l(2\sum_{k=1}^i h_{x,k} - 2 + h_{x,i+1}\r)\\
%& = & \l(l_{x,i-1}h_{x,i}^2 + n_{x,i}h_{x,i+1}^2\r) + 2(n_{x,i}h_{x,i+1}-l_{x,i-1}h_{x,i})\l(\sum_{k=1}^ih_{x,k}-1\r)\\
& = & \frac{h_{x,i}(2+h_{x,i})}{h_{x,i-1}(h_{x,i-1}+h_{x,i})} + \frac{h_{x,i+1}(h_{x,i+1}-2)}{h_{x,i}(h_{x,i-1}+h_{x,i})} + 2(h_{x,i+1}-h_{x,i})\sum_{j=1}^ih_{x,i} ~~\ge ~~ c_x(h).
\eeee
A similar argument indicates that
$$\overline{\Delta}_{y,h}w_{i,j} \ge c_y(h).$$
%Letting $c(h) = a^{-2}c_x(h)+b^{-2}c_y(h),$ it follows that $c(h)\ge 0$ by the assumptions. 
Combining these results gives that $\tilde{M}w_{i,j} \ge k(h)\boldsymbol{x},$ where $\boldsymbol{x} = (1,1,\ldots,1)^{\tT}\in\RR^{N_xN_y}.$ 

Note that $w_{i,j}$ is nonnegative on $\overline{\DD}_h$ and $\|w_{i,j}\|_\infty = 1/2.$ It follows therefore
$$\tilde{M} (\|f_{i,j}\|_\infty w_{i,j} + u_{i,j}) \ge \|f_{i,j}\|_\infty k(h)\boldsymbol{x} - f_{i,j} \ge 0,$$
based on our earlier assumptions. Hence, by the Discrete Maximum Principle, $\|f_{i,j}\|_\infty w_{i,j} + u_{i,j}$ achieves its maximum on $\partial\DD_h.$ Since $u_{i,j} = 0$ on $\partial\DD_h,$ we have 
$$u_{i,j} \le \|f_{i,j}\|_\infty w_{i,j} + u_{i,j} \le \|f_{i,j}\|_\infty \|w_{i,j}\|_\infty \le \frac{\|f_{i,j}\|_\infty}{2}$$
for $(x_i,y_j)\in\DD_h.$ Similarly, we have
$$\tilde{M}(\|f_{i,j}\|_\infty w_{i,j} - u_{i,j}) \ge \|f_{i,j}\|_\infty k(h)\boldsymbol{x} + f_{i,j} \ge 0,$$
which gives
$$-u_{i,j} \le \|f_{i,j}\|_\infty w_{i,j} - u_{i,j} \le \|f_{i,j}\|_\infty \|w_{i,j}\|_\infty \le \frac{\|f_{i,j}\|_\infty}{2}.$$
A combination of the above results gives
$$\frac{\|u_{i,j}\|_\infty}{\|f_{i,j}\|_\infty} = \frac{\|-\tilde{M}^{-1}f_{i,j}\|_\infty}{\|f_{i,j}\|_\infty} \le \frac{1}{2}.$$
It follows naturally that $\|\tilde{M}^{-1}\|_\infty \le 1/2$ since $f_{i,j}$ are arbitrary. This completes our proof.
\end{proof}

We note that our restriction on $k(h)$ is reasonable and has no negative impact on desired computations. 
In fact, $k(h)\ge 1$ can be observed if uniform or symmetric nonuniform grids are utilized. Moreover, for certain nonuniform meshes we have the following corollary.

\begin{cor}
If the nonuniform mesh $\DD_h$ is the image of a uniform mesh on $\DD$ via a smooth mapping, then we have $\|M^{-1}\|_\infty \le \sigma_{\max}/2 +\OO(\textstyle\max_{i,j}\{h_{x,i}^2,h_{y,j}^2\}),~
\textstyle\max_{i,j}\{h_{x,i}^2,h_{y,j}^2\}\rightarrow 0^+.$ 
\end{cor}

\begin{proof}
Assume that $u\in C^4(\DD)$ and $\DD_h$ is generated via smooth mappings $g^{(x)}$ and $g^{(y)}$ with
$$x_i = g^{(x)}(\omega^{(y)}_i)\quad\mbox{and}\quad y_i = g^{(y)}(\omega^{(y)}_j),$$
where $\omega^{(x)}_i = -1+ik_1,\ i=0,\ldots,N_x+1$
and
$\omega^{(y)}_j = -1+ik_2,\ j=0,\ldots,N_y+1.$
Using the Taylor expansions of the functions $g^{(x)}$ and $g^{(y)}$ in \R{kval} yields $k(h) = 1+\OO(\max\{k_1^2,k_2^2\}).$ Just as in Lemma 2.2, we may deduce that $k_1 = \OO(h_{x,i})$ and $k_2 = \OO(h_{y,j}),$ which gives $k(h) = 1+\OO(\textstyle\max_{i,j}\{h_{x,i}^2,h_{y,j}^2\}).$ Using this fact in the proof of Lemma 3.7 gives the desired result.
\end{proof}

\begin{lemma}
If $k(h)\ge 1$ then $M^{-1}[\Phi_k-I]$ is positive for sufficiently small $\tau_k\ge 0.$
\end{lemma}

\begin{proof}
Note that
$$
M^{-1}[\Phi_k-I] = \int_{t_k}^{t_{k+1}}w_k(s)\,ds,\label{w1}
$$
where
$$
w_k(s) = M^{-1}\l[\frac{1}{2}M_xE_1(s)E_2(s)E_1(s) + E_1(s)M_yE_2(s)E_1(s) + \frac{1}{2}E_1(s)E_2(s)M_xE_1(s)\r],
$$
with
$$E_1(s) = E\l(\frac{t_{k+1}-s}{2}M_x\r)\quad\mbox{and}\quad E_2(s) = E((t_{k+1}-s)M_y).$$
In a manner similar to that in \cite{Sheng2}, it follows
\bb{w2}
w_k(s) = E_1(s)E_2(s)E_1(s) + M^{-1}{\cal E}(s),
\ee
where
\bbbb
{\cal E}(s) &=& [E_1(s),M_y]E_2(s)E_1(s)+\frac{1}{2}E_1(s)[E_2(s),M_x]E_1(s)\\
%& = & \l(-\frac{1}{2}\int_{t_{k+1}}^{t_{k+1}-s}E_1(s-\xi)[A,B]E_1(\xi)\,d\xi\r)E_2(s)E_1(s)\\
%&&~~~~~~~~ + -\frac{1}{2}E_1(s)\l(\int_{t_{k+1}}^{t_{k+1}-s}E_2(s-\xi)[B,A]E_2(\xi)\,d\xi\r)E_1(s)\\
& = & \l(\frac{1}{2}\int_{t_{k+1}-s}^{t_{k+1}}E_1(s-\xi)[M_x,M_y]E_1(\xi)\,d\xi\r)E_2(s)E_1(s)\\
&& + \frac{1}{2}E_1(s)\l(\int_{t_{k+1}-s}^{t_{k+1}}E_2(s-\xi)[M_y,M_x]E_2(\xi)\,d\xi\r)E_1(s).
\eeee
Note that
$\int_{t_k}^{t_{k+1}} E_1(s)E_2(s)E_1(s)\,ds$
is positive according to Lemmas 3.1 and 3.4, so it is sufficient to show that the norm of the 
second term in \R{w2} can be made sufficiently small. To this end,
\bbbb
\l\|M^{-1}\int_{t_k}^{t_{k+1}}{\cal E}(s)\,ds\r\|_\infty & \le & \|M^{-1}\|_\infty\int_{t_k}^{t_{k+1}} \|E_1(s)\|_\infty^2\|E_2(s)\|_\infty\|[M_x,M_y]\|_\infty s\,ds\\
%& = & \frac{\tau_k^2}{2}\|M^{-1}\|_\infty\|E_1(s)\|_\infty^2\|E_2(s)\|_\infty\|[M_x,M_y]\|_\infty\\
& \le & \frac{\tau_k^2}{2}\|M^{-1}\|_\infty\|[M_x,M_y]\|_\infty
~~ \le ~~ \frac{\sigma_{\max}\tau_k^2}{4}\|[M_x,M_y]\|_\infty.
\eeee
Thus, our claim is true for sufficiently small $\tau_k>0.$
\end{proof}

\begin{rem}
Lemma 3.8 holds for all $\tau_k\ge 0$ in many realistic situations where $\|[M_x,M_y]\|_\infty=0.$ The equality is apparently
true for cases with a symmetric degeneracy term on symmetric nonuniform grids. 
Numerical experiments show that $\|[M_x,M_y]\|_\infty$ remains well-bounded.
\end{rem}

\section{Monotonicity} \clearallnum

As has been discussed in \cite{Chan2,Acker2,Acker1,Sheng4,Levine}, another fundamental property of 
solutions of \R{b1}-\R{b3} is their monotonicity. Solutions to Kawarada problems are monotonically 
increasing with respect to time on their domain of existence. Thus, it is important that this is rigorously
observed in \R{A3} as well as \R{c8},\R{c8b}.

In the following results, there are additional restrictions place on the stochastic nonlinear reaction term in order to guarantee monotonicity. Further, our results are similar to those in the non-stochastic case. We first consider \R{c8},\R{c8b}.

\begin{lemma}
If $Mv_0+g(\varepsilon,v_0) >0,$ then the solution of \R{c8},\R{c8b} is monotonically increasing with respect to 
$t\in (t_0, T_q)$ while maintaining $0\leq u<1.$
\end{lemma}

\begin{proof}
Recall \R{c8}, \R{c8b}. Let $w(t) = v'(t)$ and consider
\bbbb
&& w'(t) = [M + g_v(\varepsilon,v(t))]w(t),~~~t_0<t<T_q,\label{sdm}\\
&& w(t_0) = Mv_0+g(\varepsilon,v_0) \equiv w_0.\label{sdma}
\eeee
The solution to the above is 
$$w(t) = E\l(\int_{t_0}^t [M+g_v(\varepsilon,v(s))]ds\r)w_0,$$
where $g_v(\varepsilon,v(t))$ is the Jacobian matrix. It is also clear that $(g_v(\varepsilon,v(t)))$ is a diagonal matrix. It follows that $g_v$ is a nonnegative matrix  
since $f_u\ge 0.$ Denote
$$A(t) = \int_{t_0}^t [M+g_v(\varepsilon,v(s))]ds.$$
Hence $A_{ij} \ge 0$ for $i\neq j.$ Consequently, by Lemma 3.2, we conclude that $E(A(t))$ is positive 
for $t_0\le t<T_q.$ By Lemma 3.1, we know that $w(t)>0$ for $t_0<t<T_q$ since $w_0>0.$  
Thus, $v'(t)>0$ for $t_0<t<T_q$ and our solution is monotonically increasing with respect to time $t.$
\end{proof}

We now consider \R{A3}.

\begin{lemma}
If $Mv_0+g(\varepsilon,v_0)>0,$ then $Mv_k+g(\varepsilon,v_k)>0$ for all $k\ge 0.$ 
\end{lemma}

\begin{proof}
Let $Mv_0+g(\varepsilon,v_0)>0.$ We note for \R{A3} that
\bbbb
Mv_{k+1} &=& M\l\lbrace \Phi_kv_k + M^{-1}[\Phi_k-I]\l[\theta g(\varepsilon,v_k)+(1-\theta)g(\varepsilon,v_{k+1})\r]\r\rbrace\\
&=& \Phi_k\l[Mv_k+\l[\theta g(\varepsilon,v_k)+(1-\theta)g(\varepsilon,v_{k+1})\r]\r]-\l[\theta g(\varepsilon,v_k)
+(1-\theta)g(\varepsilon,v_{k+1})\r],
\eeee
which implies that
\bbbb
Mv_{k+1} + g(\varepsilon,v_{k+1}) & = & \Phi_k\l[Mv_k+\l[\theta g(\varepsilon,v_k)+(1-\theta)g(\varepsilon,v_{k+1})\r]\r]\\
&&-\l[\theta g(\varepsilon,v_k)+(1-\theta)g(\varepsilon,v_{k+1})\r] + g(\varepsilon,v_{k+1})\\
& > & \Phi_k\l[Mv_k+\l[\theta g(\varepsilon,v_k)+(1-\theta)g(\varepsilon,v_{k})\r]\r]\\
&&-\l[\theta g(\varepsilon,v_{k+1})+(1-\theta)g(\varepsilon,v_{k+1})\r] + g(\varepsilon,v_{k+1})\\
& = & \Phi_k[Mv_k+g(\varepsilon,v_k)]
~~ > ~~ 0
\eeee
by inductive assumptions, Lemmas 3.1 and 3.4, and the fact that $f(\varepsilon,u)$ is monotone.
\end{proof}

\begin{lemma}
Let $Mv_0+g(\varepsilon,v_0) >0$ and $\tau_k\ge 0$ be sufficiently small. Then the solution sequence of
\R{A3}, that is, $\{v_k\}_{k\ge 0},$ is monotonically increasing while maintaining $0\le v_k<1,~k\geq 0.$
\end{lemma}

\begin{proof}
From \R{A3} we have
\bbbb
v_{k+1}-v_k &=& \Phi_kv_k + M^{-1}[\Phi_k-I]g(\varepsilon,v_k)[\theta g(\varepsilon,v_k)+(1-\theta)g(\varepsilon,v_{k+1})]-v_k\\
%&=& [\Phi_k-I]v_k + M^{-1}[\Phi_k-I][\theta g(\varepsilon,v_k)+(1-\theta)g(\varepsilon,v_{k+1})]\\
%&=& M^{-1}\l\{M[\Phi_k-I]v_k + [\Phi_k-I][\theta g(\varepsilon,v_k)+(1-\theta)g(\varepsilon,v_{k+1})]\r\}\\
%&=& M^{-1}\l\{[\Phi_k-I]Mv_k + [\Phi_k-I][\theta g(\varepsilon,v_k)+(1-\theta)g(\varepsilon,v_{k+1})]\r\}\\
%&=& M^{-1}[\Phi_k-I]\{Mv_k+[\theta g(\varepsilon,v_k)+(1-\theta)g(\varepsilon,v_{k+1})]\}\\
&>& M^{-1}[\Phi_k-I]\{Mv_k+[\theta g(\varepsilon,v_k)+(1-\theta)g(\varepsilon,v_{k})]\}\\
&=& M^{-1}[\Phi_k-I]\{Mv_k+g(\varepsilon,v_k)\}
~~>~~ 0
\eeee
according to Lemmas 3.8 and 4.2. 
\end{proof}

The above results can be generalized readily to obtain the following.

\begin{thm}
Let $Mv_{k_0}+g(\varepsilon,v_{k_0})>0$ for any $k_0\geq 0$ and $\tau_k\ge 0$ be sufficiently small for $k\ge k_0.$
Then the sequence $\l\{v_k\r\}_{k\ge {k_0}}$ generated by the semi-discretized
nonuniform scheme \R{A3} increases monotonically until unity is reached or exceeded by
one or more components of the solution vector, i.e., until quenching occurs.
\end{thm}

\section{Stability} \clearallnum

Stability of different schemes for solving singular and stochastic Kawarada equations such as \R{b1}-\R{b3} 
has been a highly nontrivial and challenging concern \cite{Josh1,Josh2,Acker2,Beau1,Sheng3,Sheng4}. 
This difficulty is mostly attributed to the fact that near quenching the temporal derivative of the solution 
may blow up faster than an exponential rate \cite{Sheng5,Beau3}. Most existing investigations simply 
consider a linear stability analysis, and later the results are reconfirmed through rigorous computational 
experiments \cite{Beau1,Beau2,Sheng3,Sheng4}. Since such considerations have been proven to be
highly effective in localized computations \cite{Twi,Sheng5}, in this paper, we will first carry out a standard 
linear stability analysis for \R{A3}. 
The study will then be extend to a nonlinear case by including contributions from the singular nonlinear stochastic
reaction term. This extension of the analysis ensures that the singular and stochastic nature of the 
differential equation is incorporated near the quenching point as quenching is approached. 

\begin{lemma}
$T_{\alpha},~\alpha=x,y,$ utilized in $M$ are congruent to symmetric matrices. 
%In particular,
%$$T_\alpha = D^{-1/2}_\alpha S_\alpha D^{1/2}_\alpha\in\RR^{N_\alpha\times N_\alpha},~\alpha=x,y,$$
%where
%$$D_\alpha = \mbox{\rm diag}\l(\delta_{\alpha,1},\dots,\delta_{\alpha,N}\r),\ S_\alpha = \mbox{\rm tridiag}\l(\gamma_{\alpha,i},m_{\alpha,i},\gamma_{\alpha,i}\r)$$
%for which
%$$\delta_{\alpha,i} = \frac{h_{\alpha,i-1}+h_{\alpha,i}}{h_{\alpha,0}+h_{\alpha,1}}\quad\mbox{and}\quad \gamma_{\alpha,i} = \sqrt{n_{\alpha,i}l_{\alpha,i}}.$$
\end{lemma}

\begin{proof}
The proof is similar to that in \cite{Beau3}.
\end{proof}

\begin{lemma}
$\rho(T_\alpha)\in (-\infty,0],$ where $\rho(\cdot)$ is the spectral radius and $\alpha=x,y.$
\end{lemma}

\begin{proof}
We may consider $T_x$ since the other case is similar. By Lemma 5.1 we have that $T_x$ is congruent to a symmetric matrix, 
hence all eigenvalues of $T_x$ are real. Further, since $T_x$ is diagonally dominant with all negative diagonal elements, 
the result follows immediately. 
\end{proof}

\begin{lemma}
All eigenvalues of $M_x~\mbox{and}~M_y$ are real and negative. Further $\mu(B^{1/2}T_xB^{1/2})< 0$ and $\mu(B^{1/2}T_yB^{1/2})<0.$
\end{lemma}

%\begin{proof}
%We first consider $M_x.$ Since $M_x = BT_x = B^{1/2}B^{1/2}T_x,$ we have $B^{-1/2}M_x = B^{1/2}T_x.$ Hence
%$$B^{-1/2}M_x(B^{1/2})^{\tT} = B^{1/2}T_xB^{1/2} = B^{1/2}T_x(B^{1/2})^{\tT}$$
%because $B$ is diagonal. Further,
%$$B^{-1/2}M_xB^{1/2} = B^{-1/2}BT_xB^{1/2} = B^{1/2}T_xB^{1/2}$$
%is congruent to a symmetric matrix, since $T_x$ is congruent to a symmetric matrix. Thus, the matrices $B^{-1/2}M_xB^{1/2}$ and $T_x$ are congruent. Since the eigenvalues of $T_x$ are real and negative, the eigenvalues of $B^{-1/2}M_xB^{1/2}$ and $M_x$ must be 
%real and negative. 

%Since we have $B^{-1/2}M_xB^{1/2} = B^{1/2}T_xB^{1/2},$ then the eigenvalues of $B^{1/2}T_xB^{1/2}$ are real and negative, which gives the result. A similar argument gives the results for $T_y,~M_y,~\mbox{and}~B^{1/2}T_yB^{1/2}.$ Remark 4.1 and the fact that a matrix and its transpose share eigenvalues gives $\mu(B^{1/2}T_xB^{1/2})<0$ and $\mu(B^{1/2}T_yB^{1/2})<0.$
%\end{proof}

\begin{proof}
The proof is similar to that in \cite{Josh3}.
\end{proof}

\begin{lemma}
We have $\l\|\textstyle\prod_{j=0}^k \Phi_j\r\|_2 \le \sqrt{\kappa(B)} = \sqrt{\max\sigma_{i,j}/\min\sigma_{i,j}}$ for $\tau_j >0.$
\end{lemma}

\begin{proof}
We begin with $\Phi_{y,j} = B^{1/2}E(\tau_j B^{1/2}T_yB^{1/2})B^{-1/2},~\tau_j >0.$
For it, we have
\bb{xx1}
B^{-1/2}\Phi_{y,j}B^{1/2} ~=~ \sum_{n=0}^\infty \frac{\l(\tau_j B^{1/2}T_yB^{1/2}\r)^n}{n!}
~ = ~ E(\tau_j B^{1/2}T_yB^{1/2})~\equiv~\tilde{\Phi}_{y,j}.
\ee
By the same token,
\bb{xx2}
B^{-1/2}\Phi_{x,j}B^{1/2} ~=~ \sum_{n=0}^\infty \frac{\l((\tau_j/2) B^{1/2}T_xB^{1/2}\r)^n}{n!}
~ = ~ E((\tau_j/2) B^{1/2}T_xB^{1/2})~\equiv~\tilde{\Phi}_{x,j}.
\ee
It then follows from \R{xx1}, \R{xx2}, and Lemmas 3.5 and 5.3 immediately that
\bbbb
\l\|\prod_{j=0}^k \Phi_j\r\|_2 & = & \l\|\prod_{j=0}^k \Phi_{x,j}\Phi_{y,j}\Phi_{x,j}\r\|_2
%& = & \l\|\prod_{j=0}^k\l(B^{1/2}\tilde{\Phi}_{x,j}B^{-1/2}\r)\l(B^{1/2}\tilde{\Phi}_{y,j}B^{-1/2}\r)\l(B^{1/2}\tilde{\Phi}_{x,j}B^{-1/2}\r)\r\|_2\\
~~ = ~~ \l\|B^{1/2}\l(\prod_{j=0}^k \tilde{\Phi}_{x,j}\tilde{\Phi}_{y,j}\tilde{\Phi}_{x,j}\r)B^{-1/2}\r\|_2\\
& \le & \|B^{1/2}\|_2\|B^{-1/2}\|_2\prod_{j=0}^k E(\tau_j (\mu (B^{1/2}T_xB^{1/2})+\mu (B^{1/2}T_yB^{1/2})))
~\leq~  \sqrt{\kappa(B)},
\eeee
where $\kappa(B)$ is the condition number of the matrix $B.$
\end{proof}

By combining the above results, we arrive at the following theorem.

\begin{thm}
For $\tau_k\ge 0$ sufficiently small, the semi-adaptive nonuniform
method \R{A3} with the stochastic nonlinear reaction term frozen is unconditionally stable in the von Neumann sense
under the spectral norm. 
%that is,
%$$\|z_{k+1}\|_2 \leq C \|z_{0}\|_2,~~~k\geq 0, $$
%where $z_0=v_0-\tilde{v}_0$ is an initial error, $z_{k+1}=v_{k+1}-\tilde{v}_{k+1}$ is the
%$(k+1)$th perturbed error vector, and $C>0$ is a constant independent of $k$ and $\tau_j$ for each $0\le j\le k.$
\end{thm}

\begin{proof}
Recall \R{A3}. When the stochastic nonlinear reaction term is frozen, the error $z_{k+1}$ takes the form of
\bbb
z_{k+1} &=& \Phi_k z_k,~~~k\geq 0.\label{a}
\eee
Iterating \R{a} repeatedly to yield
$$z_{k+1}  = \prod_{j=0}^k \Phi_jz_0.$$
%\ee
The above leads to
$$\|z_{k+1}\|_2  ~\le~  \l\|\prod_{j=0}^k\Phi_j\r\|_2\|z_0\|_2 ~\le~ \sqrt{\kappa(B)}\|z_0\|_2.$$
Our proof is thus accomplished by letting $C\equiv \sqrt{\kappa(B)}$ in the above.
\end{proof}

In the circumstance when the stochastic nonlinear reaction term in \R{A3} is not frozen, we
denote $t_Q$ as the time at which numerical quenching occurs, that is, for the first time for which 
$\|v_Q\|_\infty \ge 1.$

\begin{thm}
For $\tau_k\ge 0$ sufficiently small, the semi-adaptive nonuniform 
method \R{A3} is unconditionally stable in the von Neumann sense, that is, 
for every $t_m <t_Q$ there exists a constant $C(t_m)>0$ such that
$$ \|z_{k+1}\|_2 \leq C(t_m)\|z_{0}\|_2,~~~0\le k\le m, $$
where $z_0=v_0-\tilde{v}_0$ is an initial error, $z_{k+1}=v_{k+1}-\tilde{v}_{k+1}$ is the
$(k+1)$th perturbed error vector, and $C(t_m)>0$ is a constant independent of $k$ and $\tau_j$ for each $0\le j\le k.$
\end{thm}

\begin{proof}
Recalling \R{A3}, we have
\bbbb
v_{k+1} & = & \Phi_k v_k + M^{-1}[\Phi_k-I][\theta g(\varepsilon,v_k)+(1-\theta)g(\varepsilon,v_{k+1})],~~~k\ge 0,~\theta\in[0,1].
\eeee
It follows that
\bbb
z_{k+1} & = & \Phi_k z_k + M^{-1}[\Phi_k-I][\theta(g(\varepsilon,v_k)-g(\varepsilon,\tilde{v}_k))+(1-\theta)(g(\varepsilon,v_{k+1})-g(\varepsilon,\tilde{v}_{k+1}))]\nnn\\
& = & \Phi_k z_k + M^{-1}[\Phi_k-I][\theta g_v(\varepsilon,\xi_k)z_k+(1-\theta)g_v(\varepsilon,\xi_{k+1})z_{k+1}]\nnn\\
& = & \Phi_k z_k + \l(\int_{t_k}^{t_{k+1}}w_k(s)\,ds\r)[\theta g_v(\varepsilon,\xi_k)z_k+(1-\theta)g_v(\varepsilon,\xi_{k+1})z_{k+1}],\label{stab1}
\eee
where $\xi_j = \epsilon v_j+(1-\epsilon)\tilde{v}_j,~\epsilon\in[0,1],~j=k,k+1,$ and $k\ge 0,~\theta\in[0,1].$
Recall \R{w2}. By utilizing a left-endpoint quadrature to approximate the integral in \R{stab1}, we observe that
\bbb
z_{k+1} & = & \Phi_kz_k + \l(\tau_k\Phi_k+\OO\l(\tau_k^2\r)\r)[\theta g_v(\varepsilon,\xi_k)z_k+(1-\theta)g_v(\varepsilon,
\xi_{k+1})z_{k+1}],\nnn\\
\Gamma_{1,k+1} z_{k+1} & = & \Phi_k\Gamma_{2,k}z_k,\nnn\\
z_{k+1} & = & \Gamma_{1,k+1}^{-1}\Phi_k\Gamma_{2,k}z_k,\label{stab2}
\eee
where
$\Gamma_{1,k+1} = I - \tau_k(1-\theta)\l(I+\OO\l(\tau_k\r)\r)g_v(\varepsilon,\xi_{k+1})$ and 
$\Gamma_{2,k} = I + \tau_k\theta\l(I+\OO\l(\tau_k\r)\r)g_v(\varepsilon,\xi_k).$
It follows therefore that, by employing the standard [1/0] and [0/1] Pad{\'e} approximations, respectively,
\bbbb
\Gamma_{1,k+1}^{-1} &=& E\l(\tau_k(1-\theta)\l(I+\OO\l(\tau_k\r)\r)g_v(\varepsilon,\xi_{k+1})\r)+\OO\l(\tau_k^2\r),
~~\tau_k\rightarrow 0^+,\\
\Gamma_{2,k} &=& E\l(\tau_k\theta\l(I+\OO\l(\tau_k\r)\r)g_v(\varepsilon,\xi_k)\r)+\OO\l(\tau_k^2\r), 
~~\tau_k\rightarrow 0^+.
\eeee
Apply the above to \R{stab2}. We obtain
\bbbb
z_{k+1} & = & \l\lbrace\l[E\l(\tau_k(1-\theta)\l(I+\OO\l(\tau_k\r)\r)g_v(\varepsilon,\xi_{k+1})\r)\r]\Phi_k\l[E\l(\tau_k\theta\l(I+\OO\l(\tau_k\r)\r)g_v(\varepsilon,\xi_k)\r)\r]+\OO\l(\tau_k^2\r)\r\rbrace z_k\\
%& = & \l\lbrace\prod_{j=0}^k\l[E\l(\tau_k(1-\theta)\l(I+\OO\l(\tau_j\r)\r)g_v(\varepsilon,\xi_{j+1})\r)\Phi_k E\l(\tau_k\theta\l(I+\OO\l(\tau_j\r)\r)g_v(\varepsilon,\xi_j)\r)+\OO\l(\tau_k^2\r)\r]\r\rbrace z_k\\
& = & \l\lbrace\prod_{j=0}^k E\l(\tau_j(1-\theta)\l(I+\OO\l(\tau_j\r)\r)g_v(\varepsilon,\xi_{j+1})\r)\Phi_j E\l(\tau_j\theta\l(I+\OO\l(\tau_j\r)\r)g_v(\varepsilon,\xi_j)\r)+\sum_{j=0}^k\OO\l(\tau_j^2\r)\r\rbrace z_0.
\eeee
Now, let $G(t_m) = \max_{0\le j\le m} \l\|g_v(\varepsilon,\xi_j)\r\|_2.$ Taking the norm of both sides of the above expression yields
\bbbb
\|z_{k+1}\|_2 
%&\le & \l\lbrace \l\|\prod_{j=0}^k E\l(\tau_j(1-\theta)\l(I+\OO\l(\tau_j\r)\r)g_v(\varepsilon,\xi_{j+1})\r)\Phi_j E\l(\tau_j\theta\l(I+\OO\l(\tau_j\r)\r)g_v(\varepsilon,\xi_j)\r)\r\|_2\r.\\
%&&~~~~~~~~~~~~~~~~~~~~~~~~~~~~~~~~~~~~~~~~~~~~~~~~~~~~~~~~~~~~~~~~~~~~~~~~~~~~~\l.+c_1\sum_{j=0}^k\tau_j^2\r\rbrace\|z_0\|_2\\
& \le & \l\lbrace\l\|\prod_{j=0}^k e^{\tau_j(1-\theta)(1+c_2\tau_j)G(t_m)}\Phi_je^{\tau_j\theta(1+c_3\tau_j)G(t_m)}\r\|_2+c_1\sum_{j=0}^k\tau_j^2\r\rbrace\|z_0\|_2,
\eeee
where $c_1,~c_2,~\mbox{and}~c_3$ are constants independent of $\tau_j$ and $j.$ 
Thus, for $\tau_j>0$ sufficiently small,
\bbbb
\|z_{k+1}\|_2 &\le & \l\lbrace\l\|\prod_{j=0}^k e^{2\tau_j(1-\theta)G(t_m)}\Phi_je^{2\tau_j\theta G(t_m)}\r\|_2
+\sum_{j=0}^k\tau_j\r\rbrace\|z_0\|_2\\
& \le & \l\lbrace e^{2T_qG(t_m)}\sqrt{\kappa(B)}+T_q\r\rbrace\|z_0\|_2 ~ \le ~ C(t_m)\|z_0\|_2,
\eeee
where we have used the fact that $\textstyle\sum_{j=0}^k\tau_j\le T_q.$ It follows that $G(t_m)<\infty$ 
for all $t_m<t_Q.$ This completes our proof.
\end{proof}

\section{Numerical Experiments} \clearallnum

We consider the following degenerate stochastic Kawarada model problem
\bbb
&&\sigma(x,y)u_t = \frac{1}{a^2}u_{xx}+\frac{1}{b^2}u_{yy} + \frac{\varphi(\varepsilon)}{1-u},\quad -1<x,y<1,\label{k1}\\
&&u(-1,y,t) = u(1,y,t) = 0,\quad -1\leq y\leq 1,\label{k2}\\
&&u(x,-1,t) = u(x,1,t) = 0,\quad -1\leq x\leq 1,\label{k3}\\
&&u(x,y,t_0) = u_0(x,y),\quad -1<x,y<1,\label{k4}
\eee
where $0\leq t_0< t\le T_q<\infty,$ $u_0\in C^2([-1,1]\times[-1,1]),$ and $0\le u_0\ll 1.$

Without loss of generality, we set $t_0=0,$ $a=b=2,$ and choose $u_0(x,y)= 0.001(1-\cos(2\pi x))(1-\cos(2\pi y)),\ -1\le x,y\le 1.$ 
Temporal adaptations are initiated once 
$$\max_{1\le i\leq N_x;\,1\leq j\leq N_y}u_{i,j}(t)\geq \epsilon_0=0.90,$$ 
since quenching is expected to occur in the situation \cite{Josh3,Chan1,Chan2,Beau1,Sheng5}. 
We fix the parameter $\theta = 1/2$ in \R{A3}. 
We primarily focus on the capability of the scheme \R{A3} to recover anticipated numerical solutions and, at the same time,
explore the nonlinear numerical convergence of the semi-adaptive semi-discretized scheme. Our first experiment will serve 
as a baseline approach, where nonlinear convergence results at different time levels in the absence of degenerate 
or stochastic influences will be presented. Influences of the degeneracy term on the convergence will be closely
observed in the second experiment. Standard [1/1] Pad{\'e} approximants will be utilized for matrix exponential
evaluations \cite{Iserles}. Appropriate restrictions are realized for guaranteeing the required positivity and 
monotonicity with the particular matrices involved throughout experiments \cite{Josh1,Josh2,Josh3}.

\subsection{Example 1}

\begin{center}
{\epsfig{file=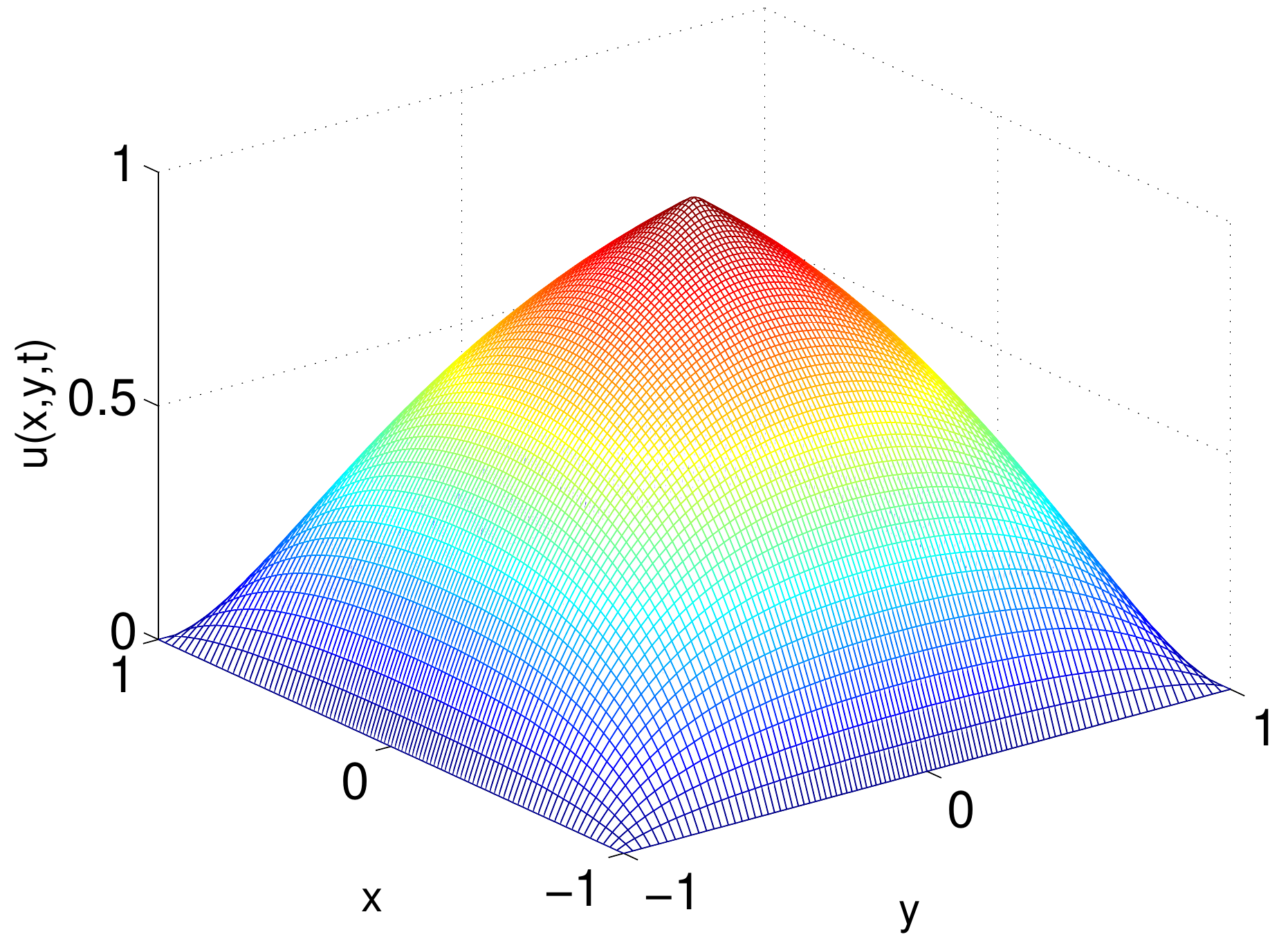,width=2.63in,height=1.68in}} %%Code: Fig_4a_4b
{\epsfig{file=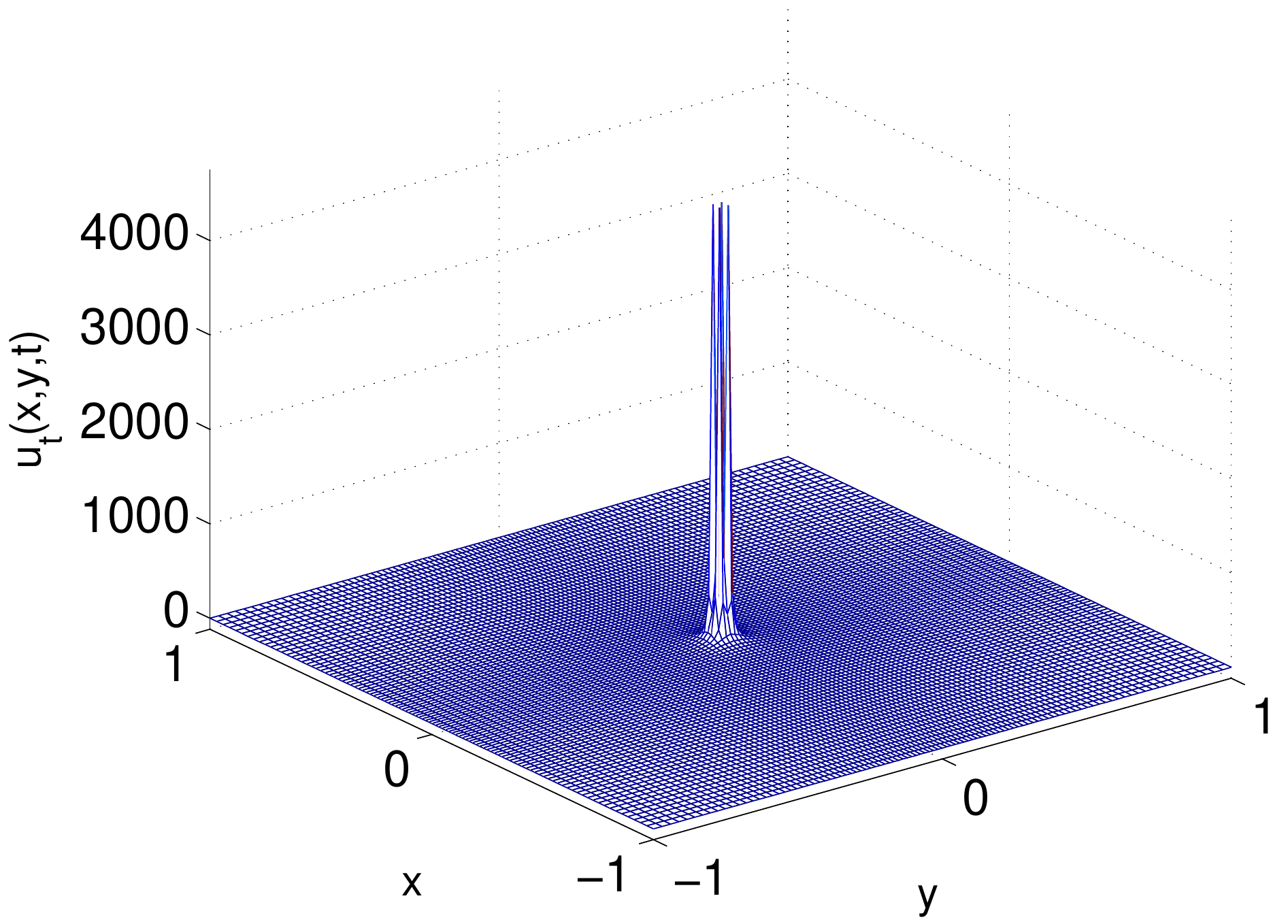,width=2.63in,height=1.68in}} %%Code: Fig_4a_4b

\parbox[t]{12.8cm}{\scriptsize{\bf Figure 1.} Solution $u$  [LEFT] and its temporal derivative $u_t$ [RIGHT]
immediately prior to quenching.}
\end{center}

To prepare our baseline results, we set $\varphi(\varepsilon) = \sigma(x,y)\equiv 1$ in \R{k1}-\R{k4}. 
Without the concerns of stochastic and degenerate impacts, we adopt a symmetric nonuniform spatial grid 
generated in a similar manner to that in \cite{Josh3}. A typical symmetric nonuniform grid can be seen in Figure 2. Adaption in the temporal direction is implemented via $u_t$ 
based arc-length monitoring functions \cite{Beau2,Sheng5,Josh3}. Quenching is observed precisely at the origin for
$T_q\approx 0.519715480937553$ which agrees well with existing results 
\cite{Josh1,Josh2,Chan3,Sheng5,Chan1,Beau1,Beau2}. While the peak value of $u$ reaches unity peacefully, 
the temporal derivative function $u_t$ grows at an explosive rate to reach $u_t^{*}\approx 4635.870128316449$ as
$t\rightarrow T_q.$ Profiles of the numerical solution and its temporal derivative function are depicted in Figure 1.

\begin{center}
{\epsfig{file=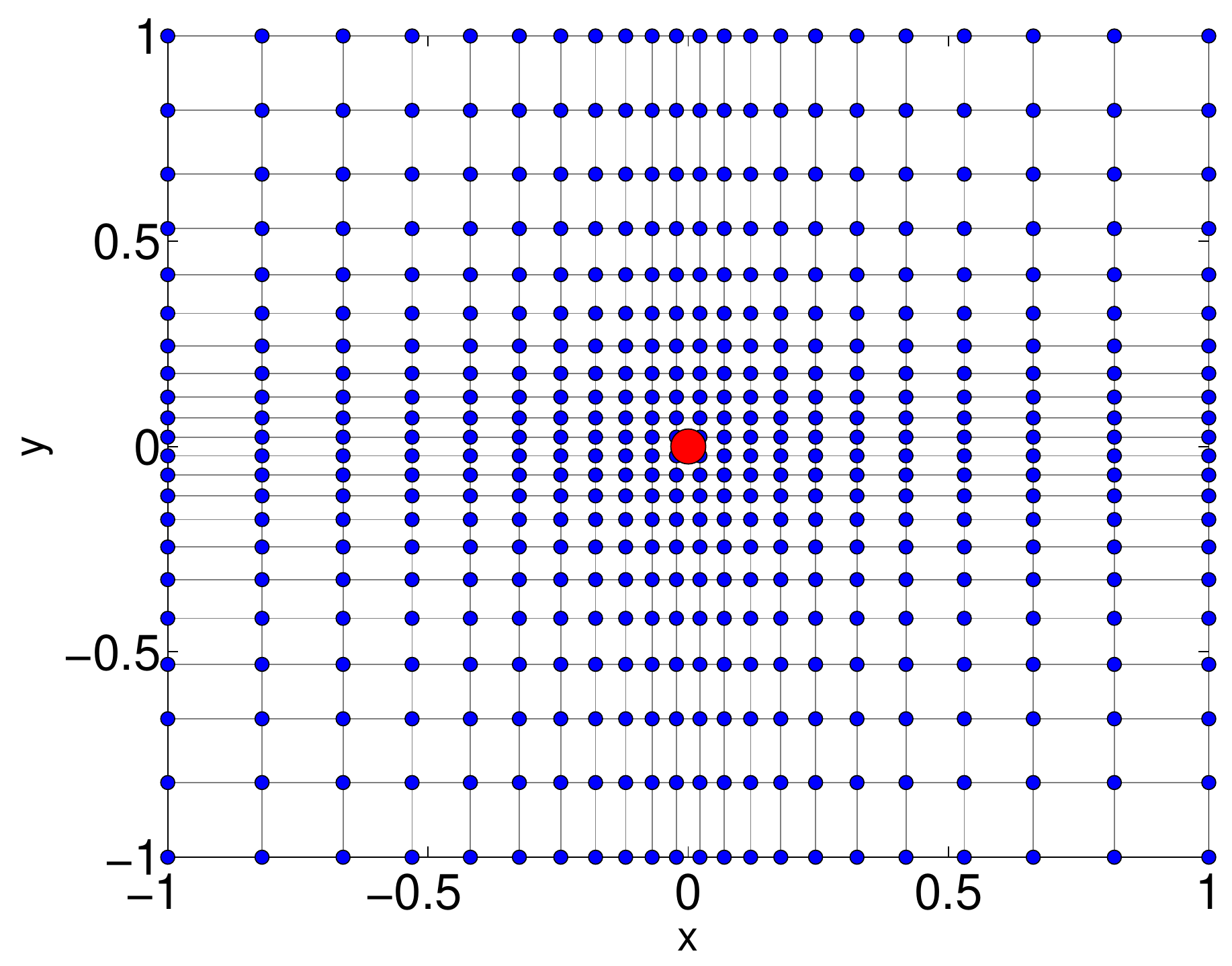,width=2.63in,height=1.68in}} %%Code: 
{\epsfig{file=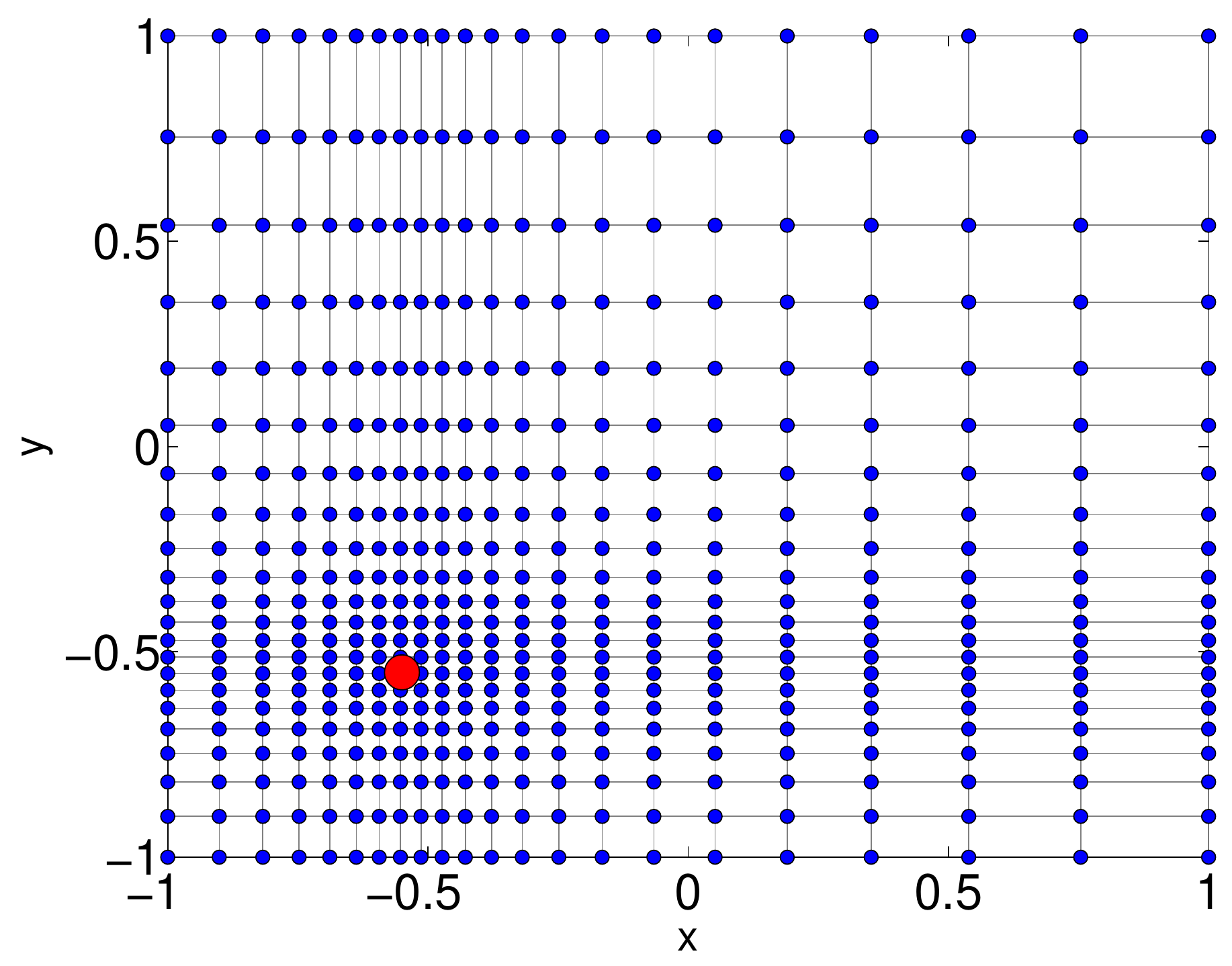,width=2.63in,height=1.68in}} %%Code: 

\parbox[t]{12.8cm}{\scriptsize{\bf Figure 2.} [LEFT] A typical symmetric nonuniform mesh used in computations for the non-degenerate case. [RIGHT] A symmetric nonuniform mesh which can be used in degenerate computations for a possible degeneracy that is concentrated at $(-1,-1).$ In both plots, the red point represents the quenching point.}
\end{center}

To examine the convergence of the numerical solution numerically, we arbitrarily choose two 
temporal points close to the quenching time $T_q:$
$$T_1 = 0.515441434291247\mbox{ and }T_2 = 0.518922378490846.$$
Let $u_h^{\tau}$ denote the numerical solution obtained by using the set of nonuniform spatial steps
represented by size $h$ and the particular adaptive temporal step $\tau.$ Suppose that $h/2$ represents
the new set of nonuniform spatial steps obtained from dividing each nonuniform mesh step in set $h$
by two. Likewise, $h/4$ is for
the new set of nonuniform spatial steps obtained from dividing each nonuniform mesh step in set $h/2$
by two. By the same token, we may define sets of temporal steps $\tau/2.~\tau/4$ based on set $\tau.$
Now, we define the following generalized Milne device for rates of point-wise convergence in space and
time, respectively:
$$p_{\rm PW}^h \approx \frac{1}{\ln 2}\ln\frac{|u_{h}^\tau - u_{h/2}^\tau|}{|u_{h/2}^\tau - u_{h/4}^\tau|},~~
q_{\rm PW}^{\tau} \approx \frac{1}{\ln 2}\ln\frac{|u_{h}^\tau - u_{h}^{\tau/2}|}{|u_{h}^{\tau/2} - u_{h}^{\tau/4}|},~~~t\in\{T_1,T_2\},$$
where these differences are defined on mesh points corresponding to $\tau,h.$ The estimates can also be 
reformulated via any vector norm. As an example of such reformulations, we include estimates for convergence rates with respect the spectral norm, denoted by $p_2^h$ and $q_2^\tau.$

%We show spectral norm based rates $p_2^h, \,q_2^{\tau}$ in Table 1 for a
%reference.

\begin{center}
\begin{tabular}{cc}
Space & Time\\
\vspace{5mm}
\begin{tabular}{l|c|c}\hline\hline
$~~~~~p_{\rm PW}^h$ & $T_1$ & $T_2$ \\ \hline
Maximum Rate    & 2.0019227 & 2.0019154 \\ \hline
Minimum Rate   & 1.8354599 & 1.8354793 \\ \hline
Median Rate   & 1.9989409 & 1.9989027 \\ \hline
Mean Rate   & 1.9983457 & 1.9980140 \\ \hline
$~~~~~p_2^h$    & 1.9976312 & 1.9955392 \\ \hline\hline
\end{tabular} %%Use data from Space_Conv.mat
& 
\begin{tabular}{l|c|c}\hline\hline
$~~~~~q_{\rm PW}^{\tau}$ & $T_1$ & $T_2$ \\ \hline
Maximum Rate  & 0.9994286 & 0.9993982 \\ \hline
Minimum Rate   & 0.9622788 & 0.8646401 \\ \hline
Median Rate   & 0.9973630 & 0.9971718 \\ \hline
Mean Rate    & 0.9944337 & 0.9913917 \\ \hline
$~~~~~q_2^{\tau}$    & 0.9831086 & 0.9536874 \\ \hline\hline
\end{tabular} %%Use data from Time_Conv.mat
\end{tabular}
\parbox[t]{12.8cm}{\scriptsize{\bf Table 1.} Detailed spatial and temporal convergence rates. 
%Spectral convergence rates in both space and time are also provided. 
The convergence
rates in space match the expectations of Lemma 2.2. Temporal convergence rates agree well with the 
Courant-Friedrichs-Lewy condition.}
\end{center}

Table 1 presents the convergence data for space and time at sampling times $T_1$ and $T_2.$ 
From Table 1 we may also notice that the mean, median, and maximal point-wise convergence rates and the spectral convergence rates
agree well in both space and time with a slight reduction when considering spectral convergence. We note that the method used to generate the nonuniform mesh for this experiment meets the criteria of Lemma 2.2, hence, we observe second-order spatial convergence. Further, these results support the Courant-Friedrichs-Lewy condition
remarkably.

\begin{center}
%{\epsfig{file=Space_Conv_Graph.eps,width=2.77in,height=1.73in}} %%Code: Fig_1a
%{\epsfig{file=Time_Conv_Graph.eps,width=2.77in,height=1.73in}}  %%Code: Fig_1b
\epsfig{file=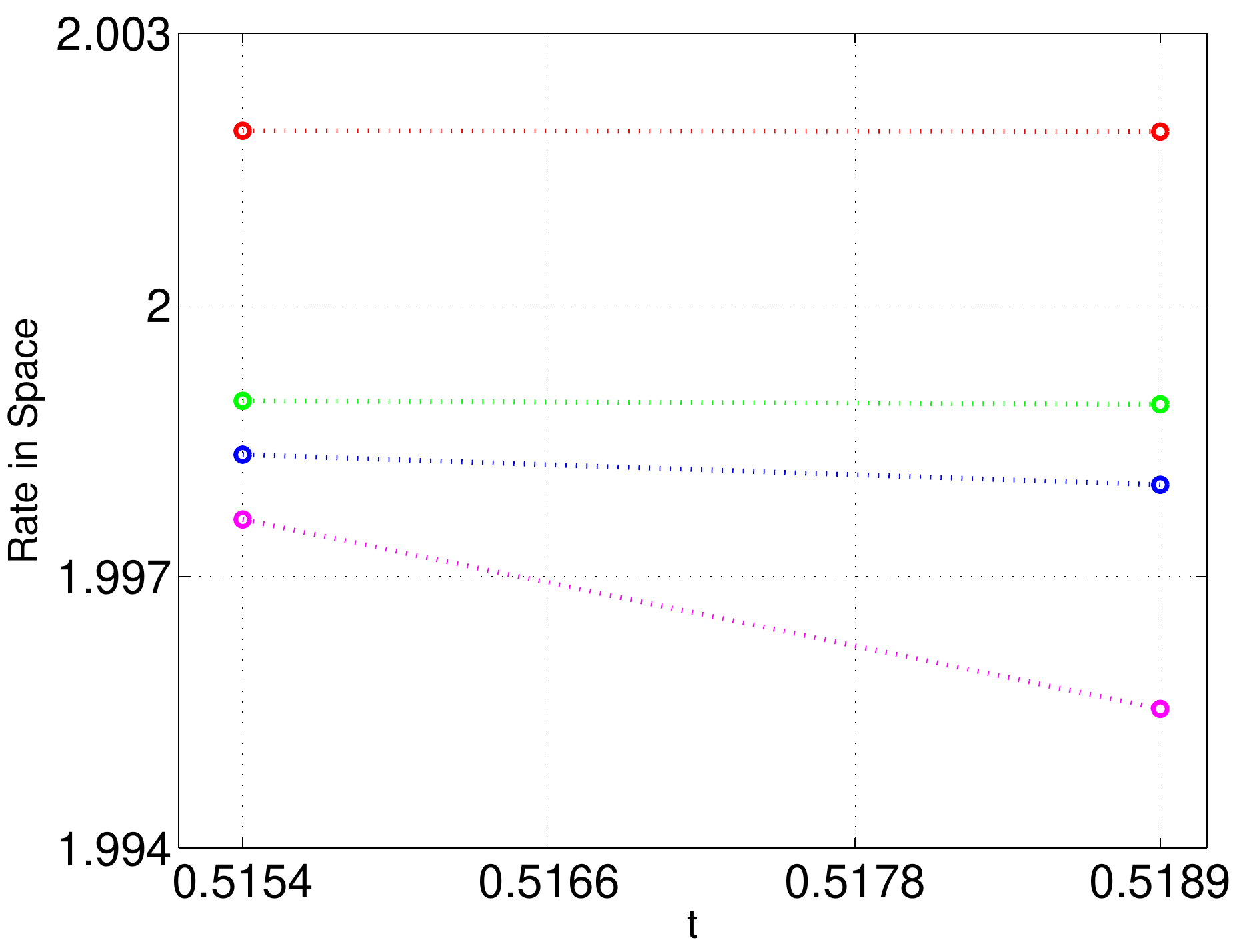,width=2.63in,height=1.68in}
\epsfig{file=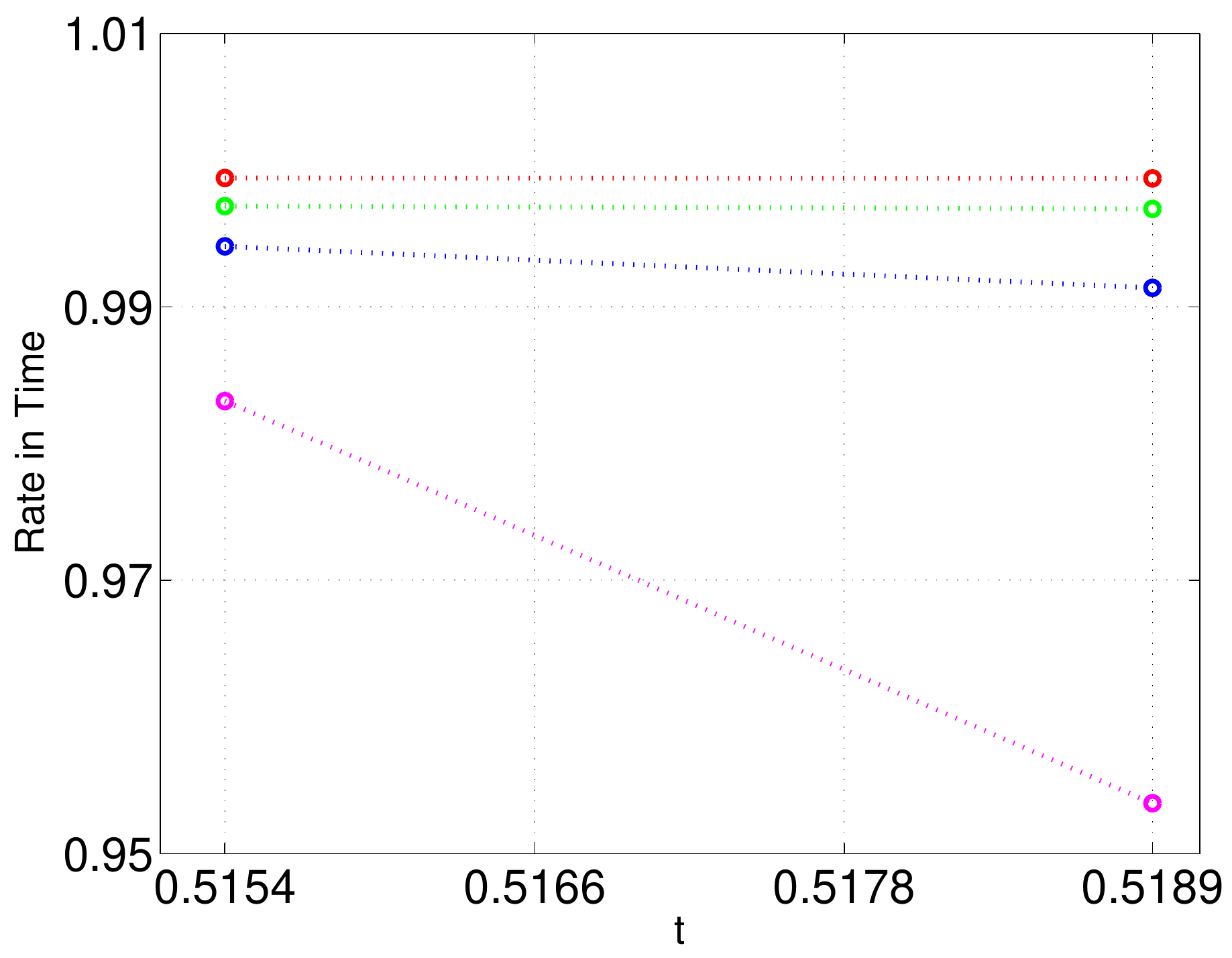,width=2.63in,height=1.68in}

\parbox[t]{12.8cm}{\scriptsize{\bf Figure 3.} Spatial convergence rate $p_{\rm PW}^h$ and $p_2^h$ [LEFT] and temporal 
convergence rate $q_{\rm PW}^{\tau}$ and $p_2^h$ [RIGHT] for the solution at times $T_1,~T_2$ are given. The
mean (blue), median (green), maximal (red), and spectral (magenta) convergence rates are included. Spatial and temporal rates are 
approximately two and one, respectively.}
\end{center}

\begin{center}
{\epsfig{file=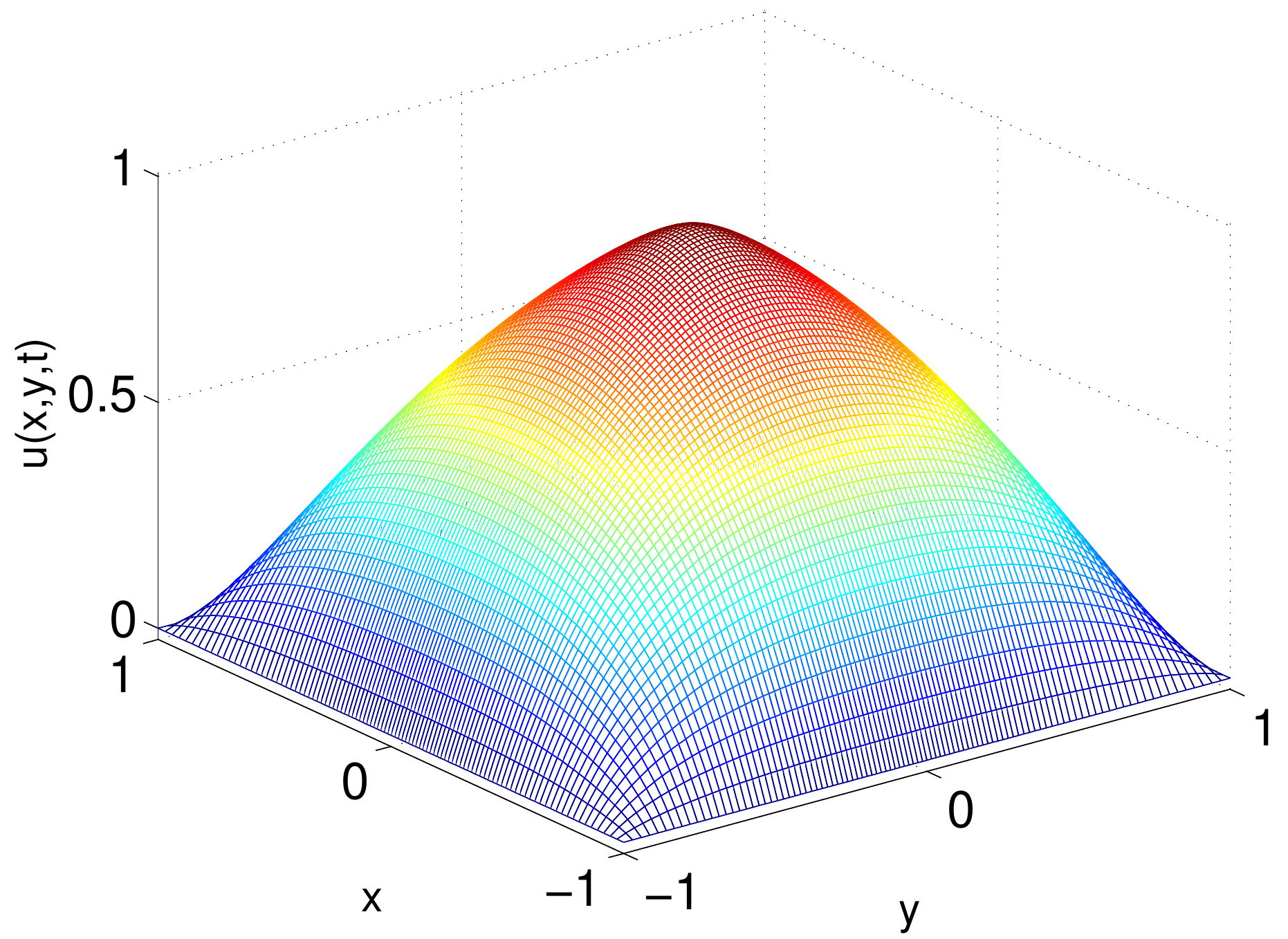,width=2.63in,height=1.68in}} %%Code: Fig_2a_2b
{\epsfig{file=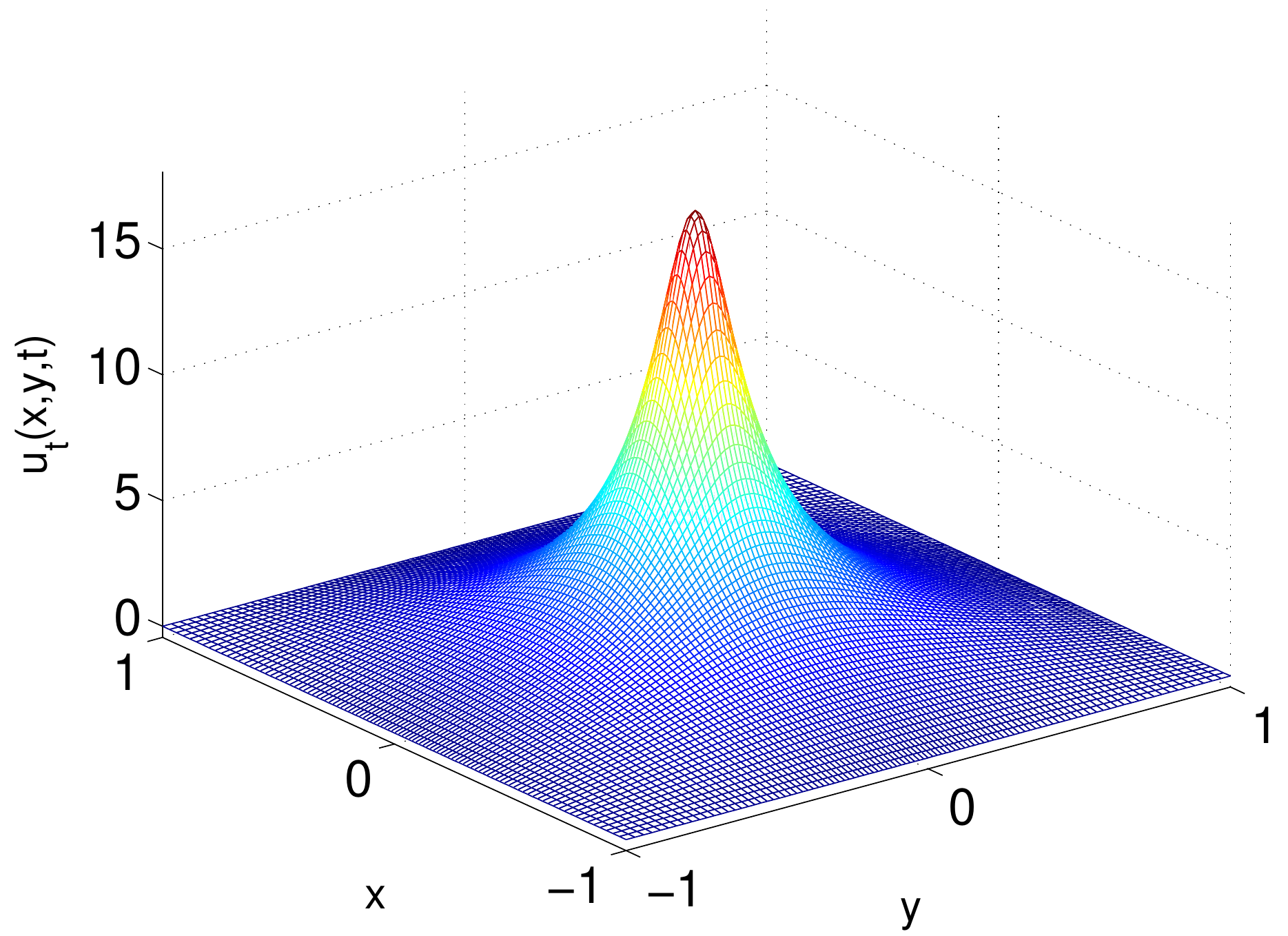,width=2.63in,height=1.68in}} %%Code: Fig_2a_2b
{\epsfig{file=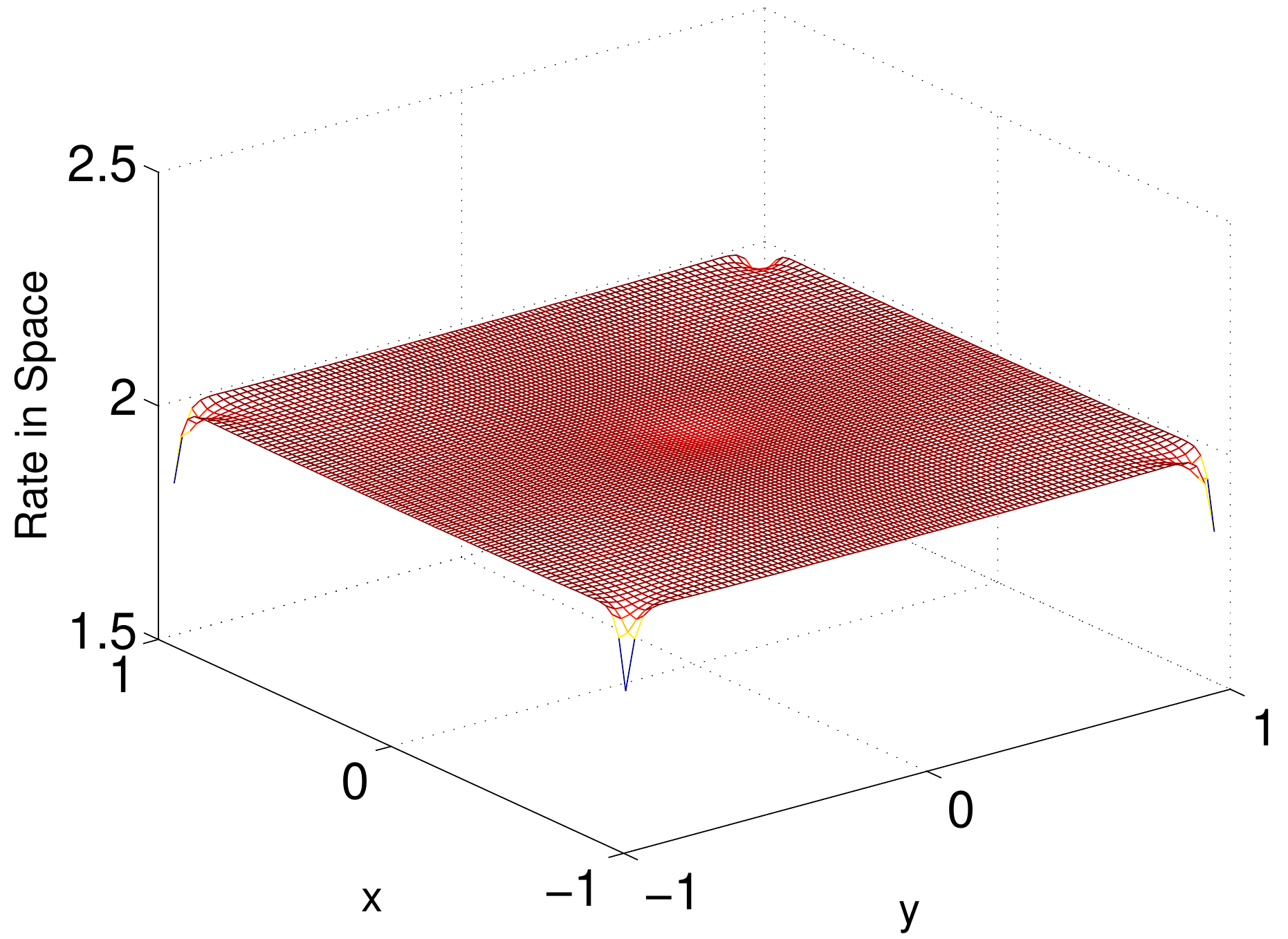,width=2.63in,height=1.68in}} %%Code: Fig_2c
%{\epsfig{file=Time_Conv_33.eps,width=2.63in,height=1.6in}} %~~ %%This is the figure with a smaller scale
{\epsfig{file=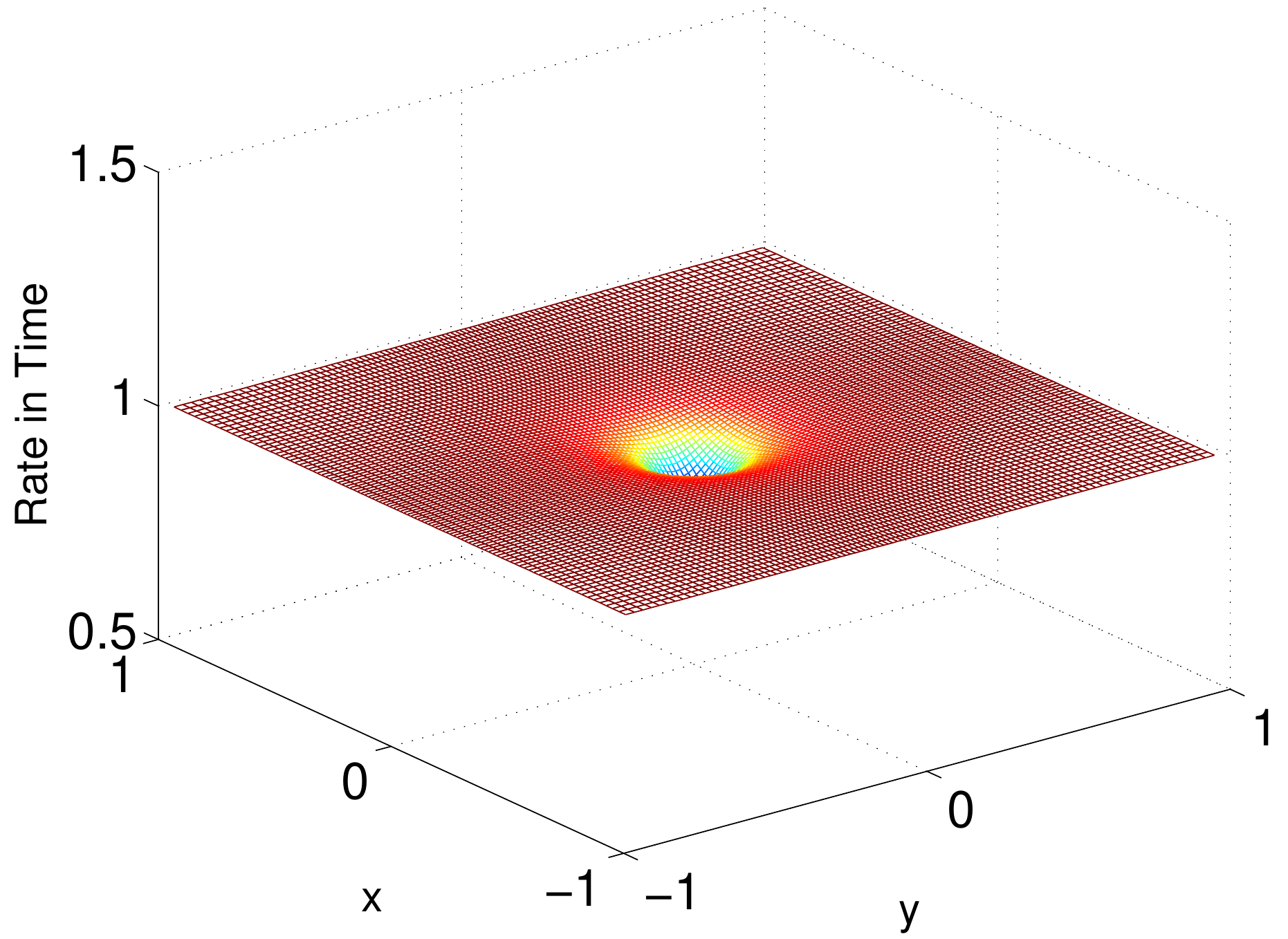,width=2.63in,height=1.68in}} %%Code: Fig_2d
\parbox[t]{12.8cm}{\scriptsize{\bf Figure 4.} 3D profiles at $T_2.$ The first row of surfaces are for the solution [LEFT] 
and its temporal derivative function [RIGHT], whereas the second row of surfaces represent the spatial convergence 
rate $p_{\rm PW}^h$ [LEFT] and temporal convergence rate $q_{\rm PW}^{\tau}$ [RIGHT].}
\end{center}

Figure 3 depicts the dynamic profile of convergence rates between $T_1$ and $T_2.$ The lines indicate movements of
the maximal, median, and mean point-wise convergence rates as well as the spectral convergence rates in both space and time. When considering the
spatial convergence, we observe that the curves are very close to two, and when considering the temporal convergence, 
we notice that the curves are grouped very close to one. Minimal spatial and temporal convergence rates are recorded in Table 1. The minimum spatial convergence rates most likely
due to the smoothness
of the solution about corners of the spatial domain \cite{KinNA}. The minimal temporal convergence rates are most likely
due to the development of quenching singularity, which indicated by rate drops from $T_1$ to $T_2.$
The quenching singularity also affects the spatial convergence, but less violently as the gradient of the solution is bounded. 
These expectations follow from a generalization of the results in \cite{Levine2}. 

Further, Figure 4 depicts the solution, its temporal derivative, and convergence rates at $T_2.$ We note that the peaks 
of the solution and the derivative are concentrated at the origin, which is the quenching location in this particular 
problem \cite{Josh3,Sheng5,Chan2,Kawa,Levine2}. In the bottom row of the figure, rates of point-wise spatial and temporal convergences 
are plotted as surfaces. It is easy to observe that both rates decrease near the origin. This, once again, is due to the 
reaction term becoming unbounded near quenching, or combustion. The singularity has an explosive impact to 
the overall convergence rates. It is also observed that that there are decreased spatial convergence rates near corners 
of the domain, which has been predicted by several authors  \cite{KinNA,Sheng6,Chan4}.

Despite drops in convergence rates near the quenching point, the mean convergence rate of the scheme is 
well-maintained near quenching at time $T_2.$ 
At this time, the mean spatial and temporal convergence rates are approximately 1.998013984025239 and 0.991391747823586, respectively.
%At this time, the mean spatial rate is approximately 1.998013984025239 
%and the mean temporal rate is approximately 0.991391747823586. 
The spectral convergence rates in space and time are 1.995539227366137 and 0.953687356004341, respectively. We note that the spatial convergence matches the expectations of Lemma 2.2 but the temporal convergence rate is lower than the expected rate from Lemma 2.1. 
%We further note that both convergence rates 
%are better than those obtained in the analysis. This is probably due to the use of a particular example with 
%symmetric nonuniform grids. 
However, the temporal convergence rate is appropriately proportional to the spatial rate, hence, the reduced temporal convergence is a result of the Courant-Friedrichs-Lewy condition.

\begin{center}
%{\epsfig{file=ut_Convergence_2.eps,width=2.63in,height=1.6in}} %~~%%Old figure
%{\epsfig{file=ut_Convergence_3.eps,width=2.63in,height=1.6in}} %~~%%Old figure
{\epsfig{file=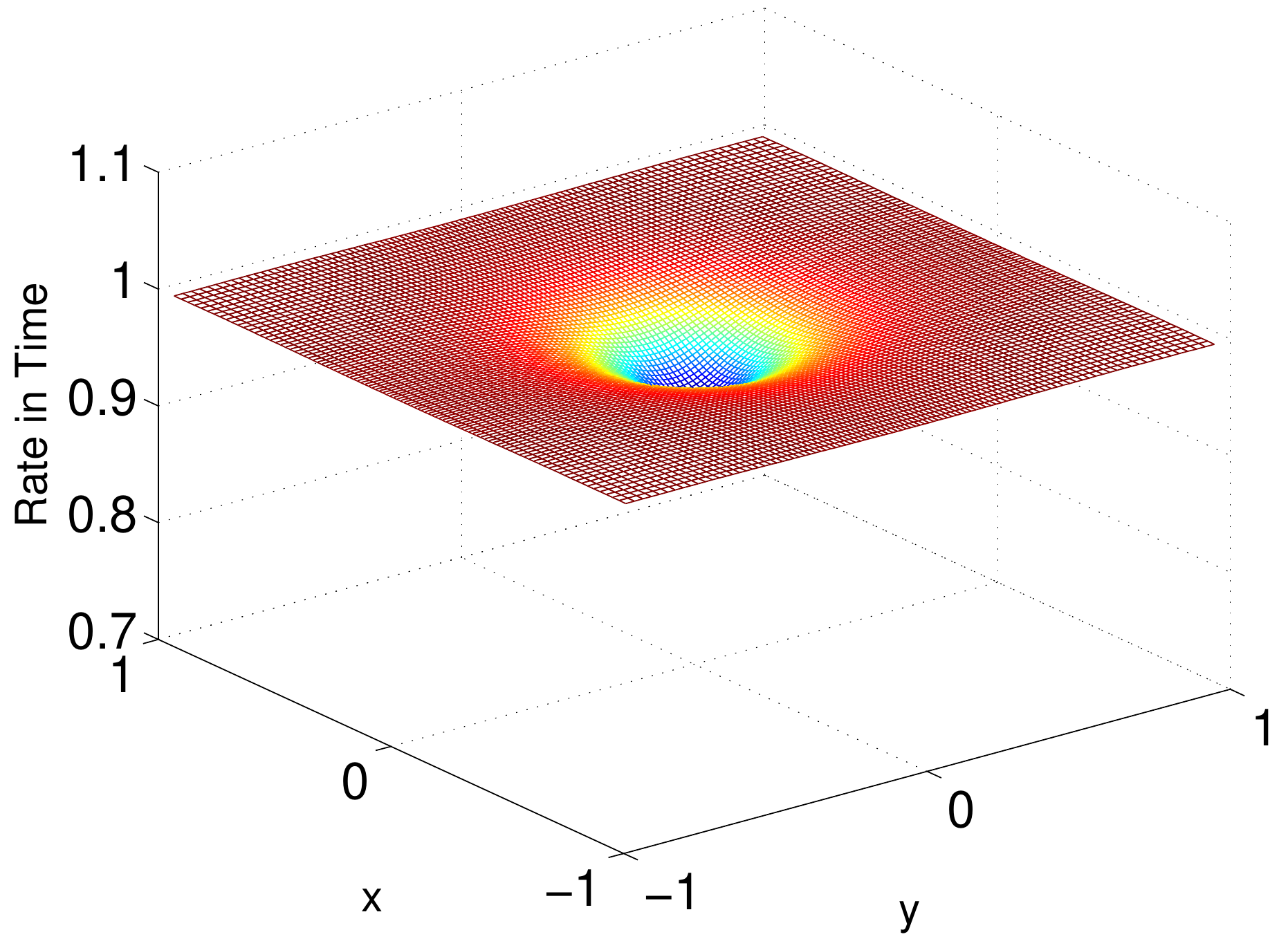,width=2.63in,height=1.68in}} %%Code: Fig_3a_3b
{\epsfig{file=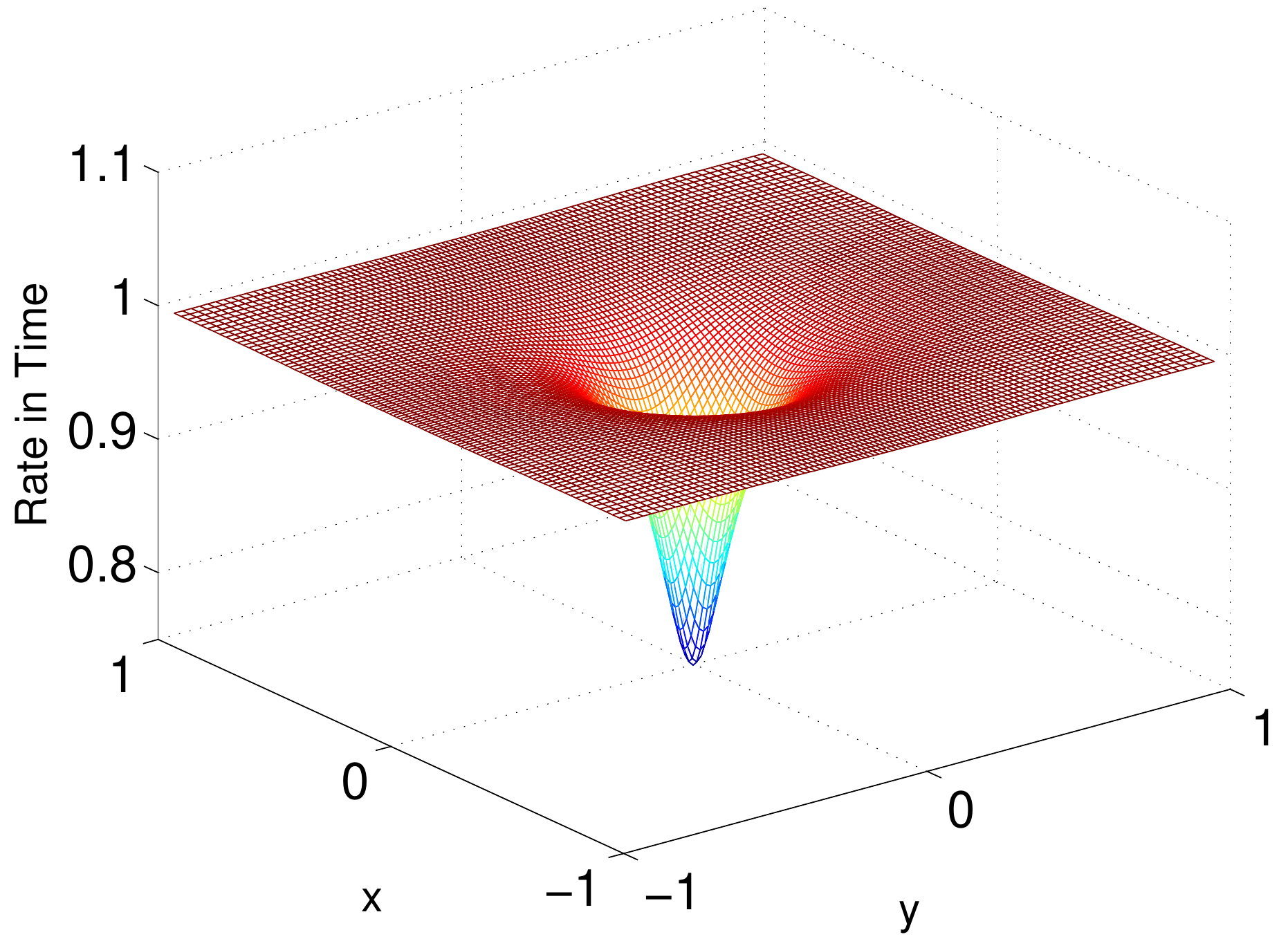,width=2.63in,height=1.68in}} %%Code: Fig_3a_3b
\parbox[t]{12.8cm}{\scriptsize{\bf Figure 5.} Surface profiles representing convergence rates,
$p_{\rm PW}^h,~q_{\rm PW}^{\tau},$ of the temporal derivative of the solution, $u_t,$ at 
$T_1$ [LEFT] and $T_2$ [RIGHT] . The minimum in the latter case is approximately 0.751476025769804.}
\end{center}

We would like to also note that implementing uniform grids or considerations of less exotic initial data will not 
drastically change the convergence rates nor the overall shape of the convergence surface plots. 
Furthermore, similar conclusions are confirmed in our experiments when nonsymmetric grids are utilized,
though these results are admitted for brevity. The effects of nonuniform grids, that are not generated in a 
smooth manner, on convergence rates have yet to be explored in relation to quenching-combustion problems. 
On the other hand, to further raise the accuracy of schemes, say, to obtain a second-order temporal accuracy, 
we would need to implement spatial discretizations to the fourth order. This could be fulfilled via a nine-point scheme or a 
compact algorithm \cite{Iserles}. 
%We do note that there is a slight reduction of convergence rate as the quenching singularity becomes more significant, but the spatial and temporal convergence rates are still approximately two and one, respectively, as the influence increases. Moreover, it is observed that quenching has a much larger impact on temporal convergence than spatial convergence. 
We also observe that, while the emergence of the quenching singularity reduces the spatial convergence rate 
significantly near the origin, the reduction is less severe than the reduction in temporal convergence. This is probably 
due to the fact that $u$ and its spatial derivatives are well-behaved in our quenching computations \cite{Levine2,Sheng8}.

The temporal derivative of the solution of \R{k1}-\R{k4} becomes unbounded as $t\rightarrow T_q^-$ \cite{Sheng8}. 
As a reference, in Figure 5 we explore the temporal convergence rate of $u_t$ at times $T_1$ and $T_2,$ 
respectively, in order to examine whether the explosive rate of change can significantly affect the convergence rate. 
The shapes of the surfaces in Figure 5 are no doubt similar to that in Figure 4. We find that the convergence 
rates drop significantly near the origin. However, the mean temporal convergence rate is still approximately first-order. 
We observe the average point-wise convergence rate is approximately 0.991674641553278 at $T_1$ and 
0.986792663593675 at $T_2.$ Further, the spectral convergence rates of $u_t$ are 0.948154262871667 and 0.819817700033525 at times $T_1$ and $T_2,$ respectively. This novel exploration of convergence issues provides profound insights into the full 
nature of singular nonlinear problems such as \R{b1}-\R{b3}.

\subsection{Example 2}

In this experiment we explore impacts of a degeneracy and a mild stochastic influence on both the solution and rates of convergence.
In the following, we set $\varphi(\varepsilon)=\varepsilon$ for $\varepsilon\in[0.98,1.02],$ and $\sigma(x,y) = \textstyle\sqrt{(x+1)^2+(y+1)^2},\ -1\le x,y\le 1.$ 
%in order to observe impacts of the degeneracy on both the solution and rates of convergence. 
We note that this choice 
for $\sigma(x,y)$ forces the degeneracy to occur at the spatial point $(x_d,y_d)=(-1,-1)$ \cite{Sheng3,Beau3}. By choosing such a degeneracy 
function, we will be better able to observe its fast spreading effects on the solution. Our nonuniform grids are generated 
in the consideration that the quenching location will be shifted, but not too far away, from the origin \cite{Sheng5}. See Figure 2 for a possible example of how to shift the mesh appropriately.
The convergence rates are computed in a similar fashion as before. 

\begin{center}
{\epsfig{file=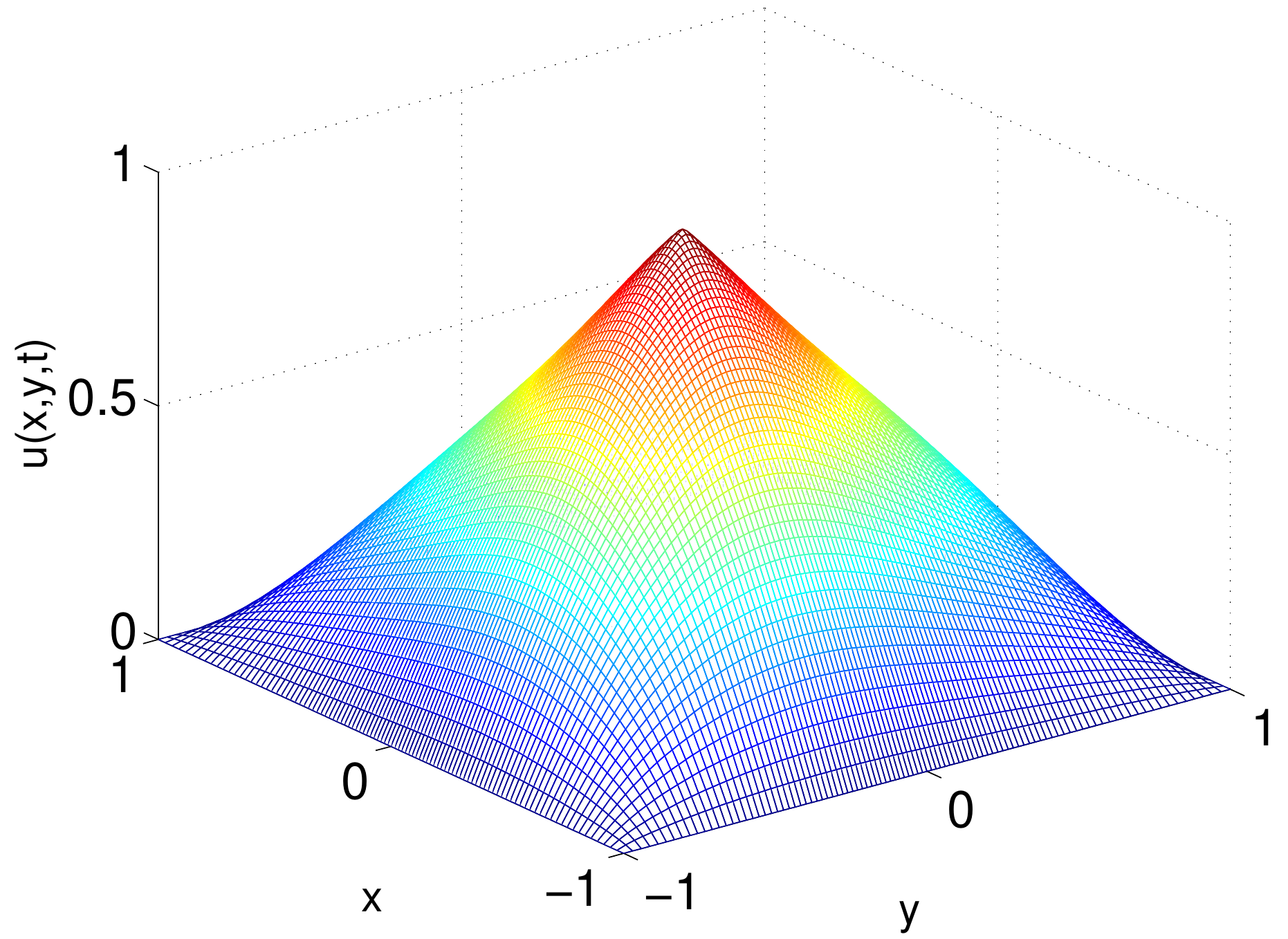,width=2.63in,height=1.68in}} %%Code: Fig_8a_8b
{\epsfig{file=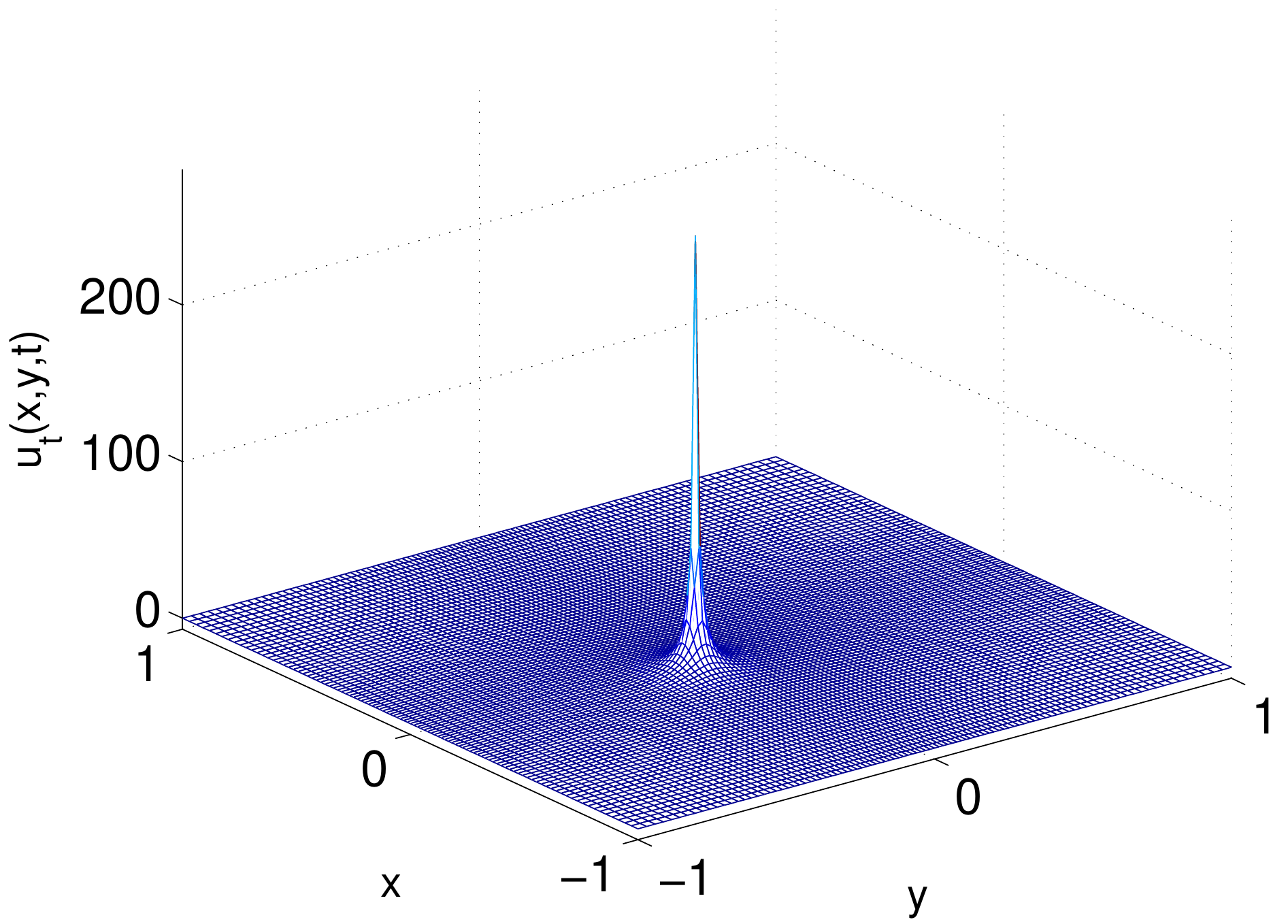,width=2.63in,height=1.68in}} %%Code: Fig_8a_8b
\parbox[t]{12.8cm}{\scriptsize{\bf Figure 6.} Solution $u$ [LEFT] and its temporal derivative function $u_t$ [RIGHT]
as $t\rightarrow T_q.$ We have the maximum $\textstyle_{-1\le x,y\le 1} u_t(x,y,t) =
u_t(x_q,y_q,t) \approx 279.5375137783287$ as $t\rightarrow T_q.$}
\end{center}

An elevated quenching time, $T_q \approx 0.836473787737391,$ probably primarily due to the degeneracy \cite{Chan4}, 
is observed. However, the vibrant stochastic feature of the reaction term may also increase quenching time \cite{Josh3}. Figure 6 depicts the solution and its derivative immediately prior to quenching. Just as we have
anticipated, the quenching location is slightly shifted from the origin to $(x_q,y_q)\approx 
(-0.151533166458073,-0.151533166458073).$ The peak value of $u_t,$ which is concentrated around the quenching 
location, reaches 279.5375137783287 prior to quenching. 

Due to our slightly delayed ignition-quenching time $T_q,$ we consider the following two arbitrary convergence testing times
$$T_1 = 0.827203763602093\mbox{ and }T_2 = 0.834206506292535.$$
A major change in spatial convergence is expected since the degeneracy $\sigma(x,y)$ is a spatial function 
creating strong singularities on the boundary of the spatial domain.
In particular, a modification of the proof of Lemma 2.2 will show that the degeneracy limits the effects of generating a nonuniform mesh via a smooth mapping. It seems that the stronger the singularity is, the more the order of spatial accuracy and convergence is reduced.

\begin{center}
\begin{tabular}{cc}
Space & Time\\
\vspace{5mm}
\begin{tabular}{l|c|c}\hline\hline
$~~~~~p_{\rm PW}^h$ & $T_1$ & $T_2$ \\ \hline
Maximum Rate  & 3.8371304 & 3.7599219 \\ \hline
Minimum Rate   & -1.2521872 & -1.2048243 \\ \hline
Median Rate  & 1.1193987 & 1.1231935 \\ \hline
Mean Rate  & 1.1156642 & 1.1258041 \\ \hline
$~~~~~p_2^h$    & 1.2147392 & 1.3080651 \\ \hline\hline
\end{tabular} 
& 
\begin{tabular}{l|c|c}\hline\hline
$~~~~~q_{\rm PW}^{\tau}$ & $T_1$ & $T_2$ \\ \hline
Maximum Rate   & 0.9999996 & 0.9999993 \\ \hline
Minimum Rate  & 0.9989553 & 0.9956911 \\ \hline
Median Rate     & 0.9999220 & 0.9999153 \\ \hline
Mean Rate        & 0.9998710 & 0.9998046 \\ \hline
$~~~~~q_2^{\tau}$    & 0.9995514 & 0.9985600 \\ \hline\hline
\end{tabular}
\end{tabular}
\parbox[t]{12.8cm}{\scriptsize{\bf Table 2.} Detailed spatial and temporal convergence rates. The experimental 
results suggest clearly a slightly above one spatial convergence rate, and slightly below one temporal rate.}
\end{center}

Table 2 is devoted to convergence information in space and time at sampling times $T_1$ and $T_2.$ 
From Table 2 we may also notice that while the mean and median convergence rates strongly suggest
a slightly higher than one rate of convergence in space and slightly lower than one rate of convergence
in time, minimal and maximal rate estimates differ significantly. Again, the point-wise rates of convergence and temporal spectral rate of convergence
decrease as $t\rightarrow T_q^-,$ though reductions are less significant probably also due to a 
stronger restriction on temporal step sizes to maintain the solution positivity when $\sigma$ 
is nontrivial \cite{Josh3,Josh2,Josh1,Sheng5,Beau1,Beau2,Beau3}. 
The spatial spectral convergence rate actually increases in time, but this is due to the solution becoming more smooth near the corners of the domain as the solution is iterated.
The phenomena are well
documented by Figure 7 where spatial and temporal convergence rates of the numerical solution 
under the influences of a non-constant $\sigma(x,y)$ and a nonlinear stochastic reaction term are exhibited.

\begin{center}
%{\epsfig{file=Deg_Space_Convergence_Graph1.eps,width=2.63in,height=1.6in}} %~~
%{\epsfig{file=Deg_Space_Conv_Graph.eps,width=2.77in,height=1.73in}}  %%Code: Fig_5a
%{\epsfig{file=Deg_Time_Convergence_Graph1.eps,width=2.63in,height=1.6in}} %~~
%{\epsfig{file=Deg_Time_Conv_Graph.eps,width=2.77in,height=1.73in}} 
\epsfig{file=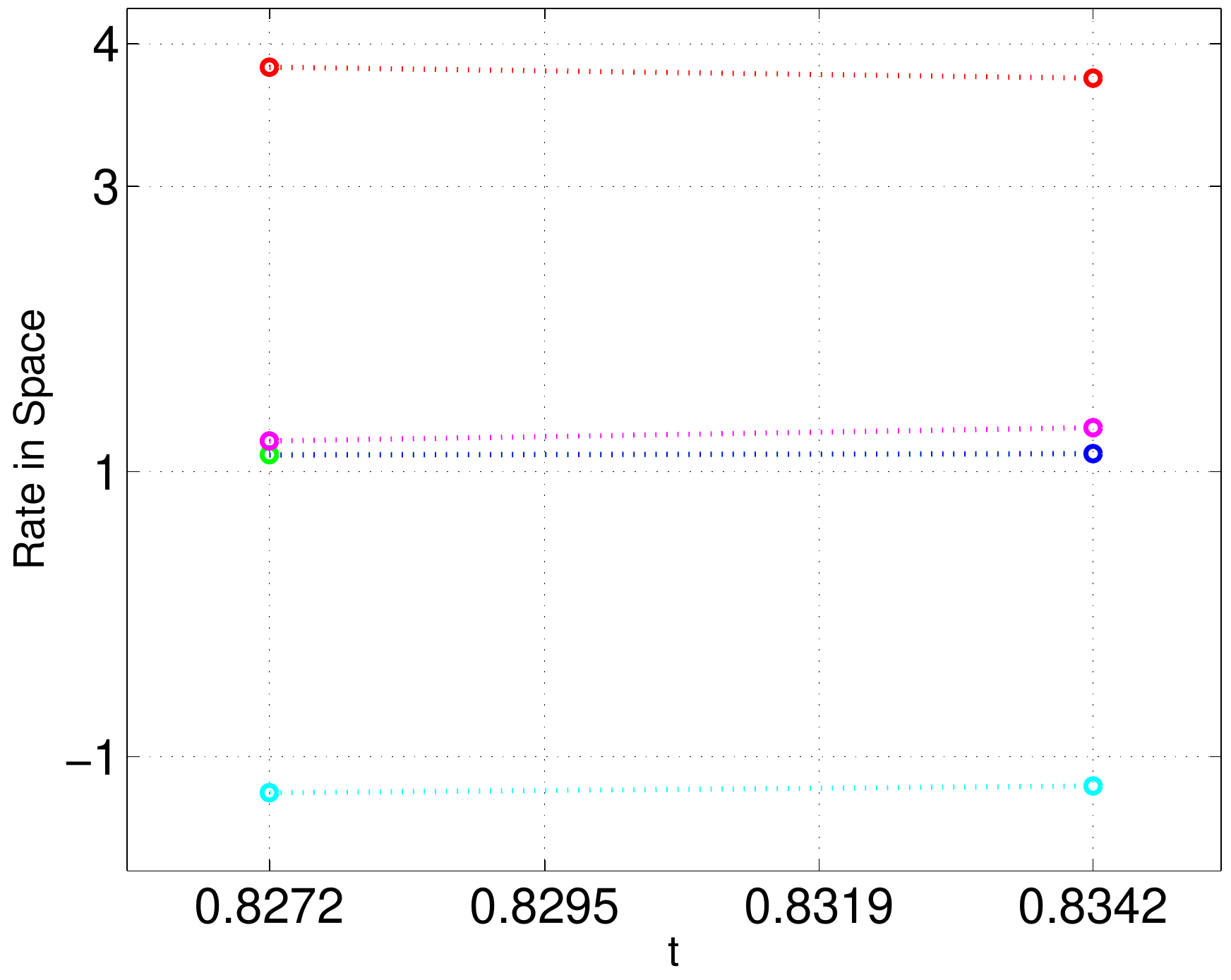,width=2.63in,height=1.68in}
\epsfig{file=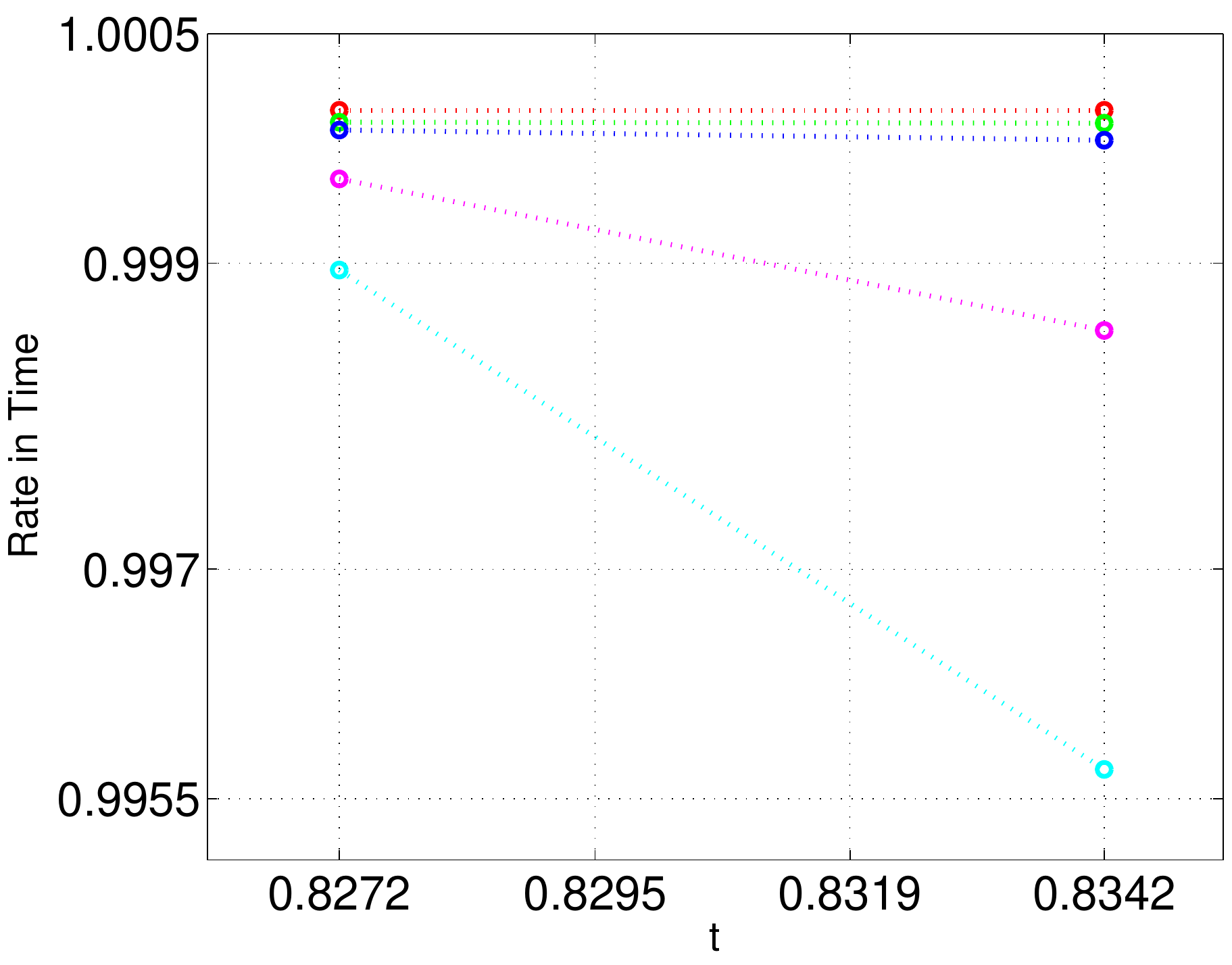,width=2.63in,height=1.68in}
\parbox[t]{12.8cm}{\scriptsize{\bf Figure 7.} Spatial convergence rates $p_{\rm PW}^h$ and $p_2^h$ [LEFT] and temporal convergence rates $q_{\rm PW}^h$ and $q_2^h$ [RIGHT] 
for the solution at times $T_1$ and $T_2.$ The mean (blue), median (green), maximal (red), minimal (cyan), and spectral (magenta) rates are included. The mean and
median curves almost overlap in the first picture. The negative 
minimal rate of convergence in space is due to spurious artifacts of the degeneracy.}
\end{center}

\begin{center}
{\epsfig{file=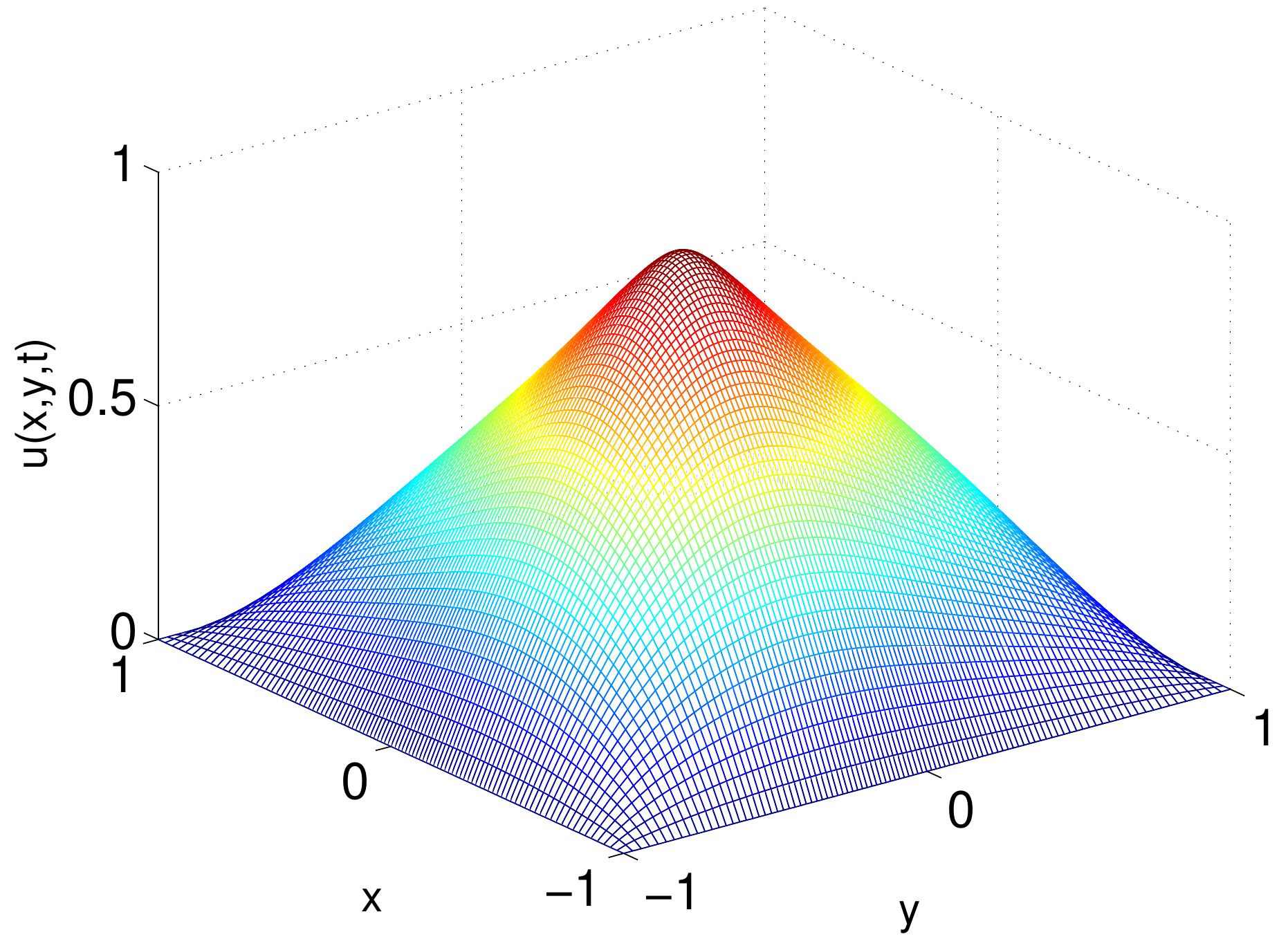,width=2.63in,height=1.68in}}  %%Code: Fig_6a_6b
{\epsfig{file=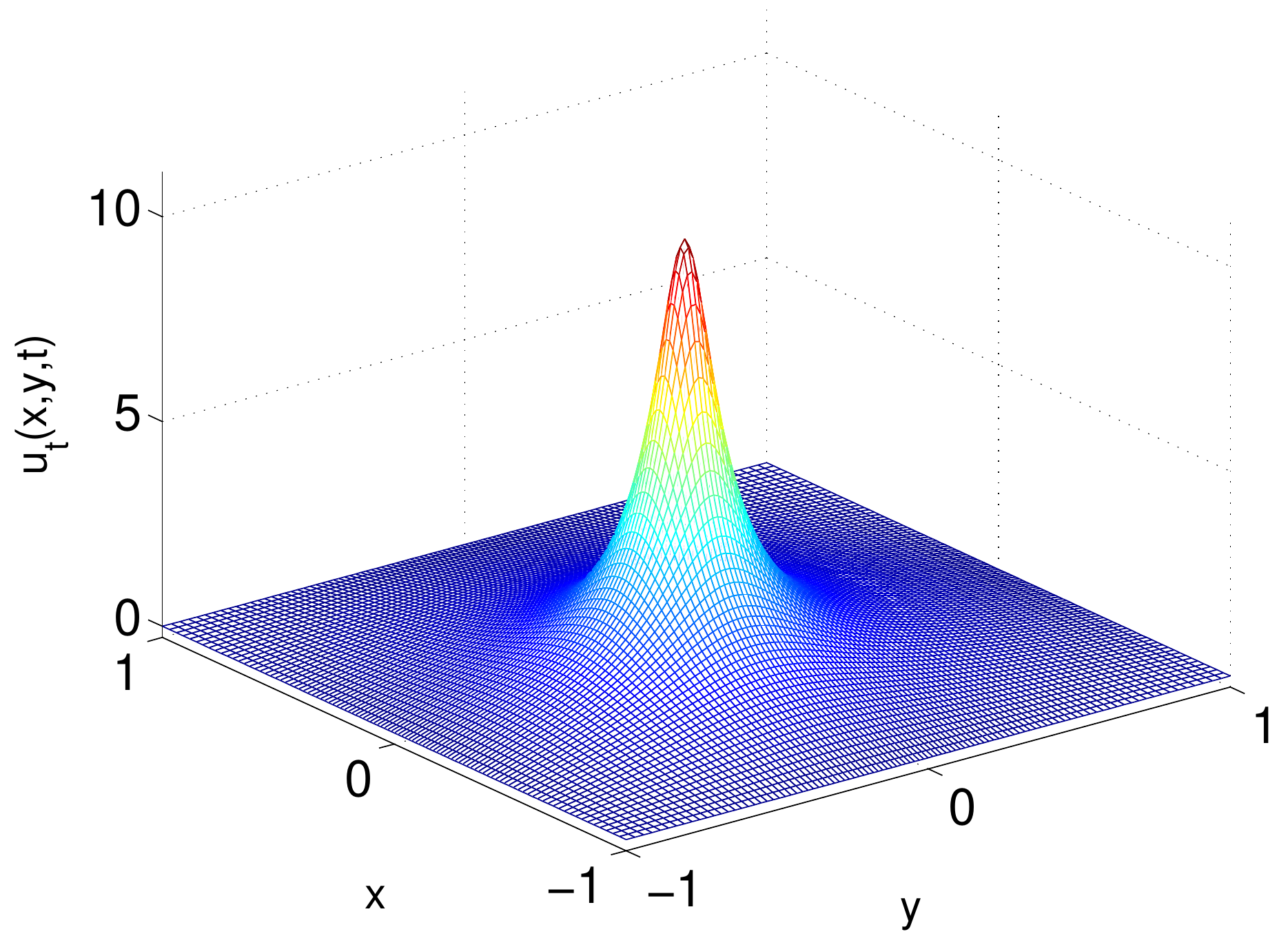,width=2.63in,height=1.68in}}  %%Code: Fig_6a_6b
{\epsfig{file=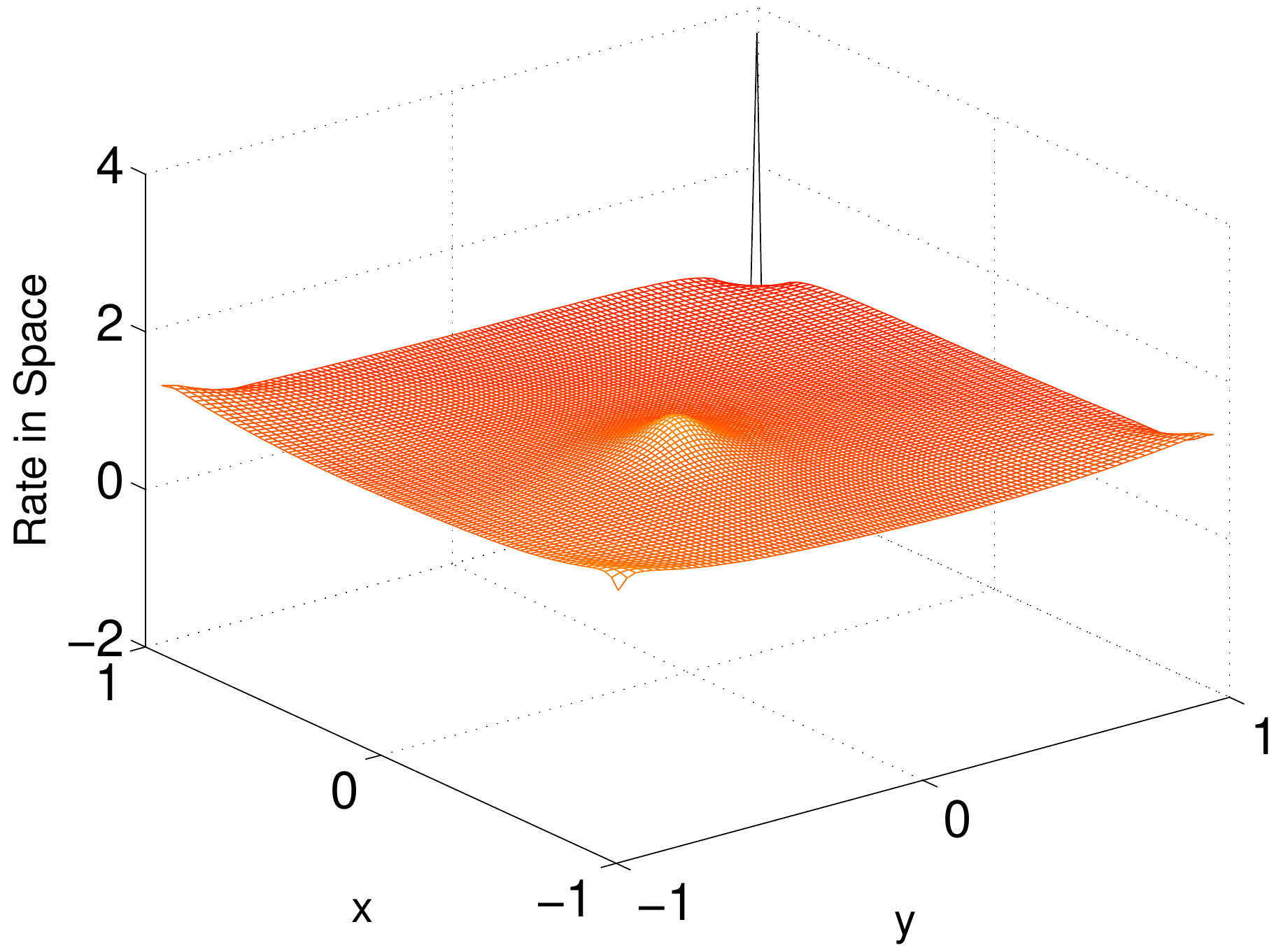,width=2.63in,height=1.68in}}   %%Code: Fig_6c
{\epsfig{file=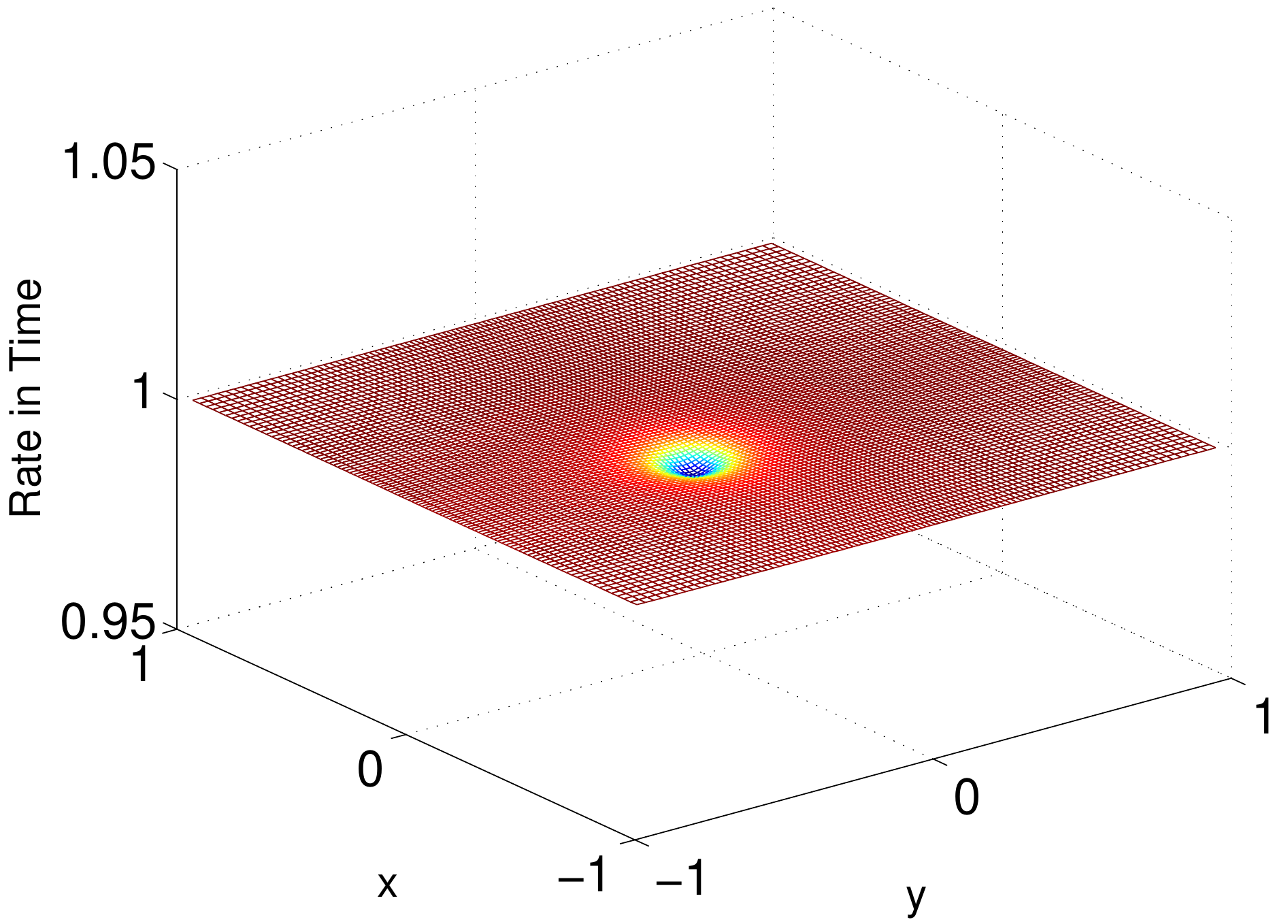,width=2.63in,height=1.68in}}      %%Code: Fig_6d
\parbox[t]{12.8cm}{\scriptsize{\bf Figure 8.} Surface profiles at time $T_2.$ The first row is for
the solution [LEFT] and its temporal derivative [RIGHT], whereas the second row is for 
the rates of convergences in space [LEFT] and in time [RIGHT].}
\end{center}

Figure 8 illustrates the solution of our stochastic degenerate problem, as well as its derivative, at time $T_2.$ 
Surface plots of the spatial and temporal convergence are also presented. We note that there is a visible shift of
the quenching location, as well as the position of the peaks of $u,~u_t,$ to $(x_q,y_q).$ This shift agrees well 
with existing predictions \cite{Sheng5}. Slight dips at the quenching point can also be observed in Figure 8.
Needless to say, the point-wise convergence rate surfaces provide us with further insights into the maximal and minimal
rates noted in Table 2 and Figure 7. From the surface plot, we observe that spurious issues occur at the corner opposite 
of the degeneracy in the spatial domain. This phenomenon can be explained by the fact that the degeneracy simplifies 
the problem at $(-1,-1)$ by reducing the parabolic problem computationally to an elliptic problem. However, the 
degeneracy still creates disturbances, which are observed in the extreme convergence rates near $(1,1).$ 

\begin{center}
%{\epsfig{file=ut_Convergence_2_degeneracy_time.eps,width=2.63in,height=1.6in}} %~~
%{\epsfig{file=ut_Convergence_3_degeneracy_time.eps,width=2.63in,height=1.6in}} %~~
{\epsfig{file=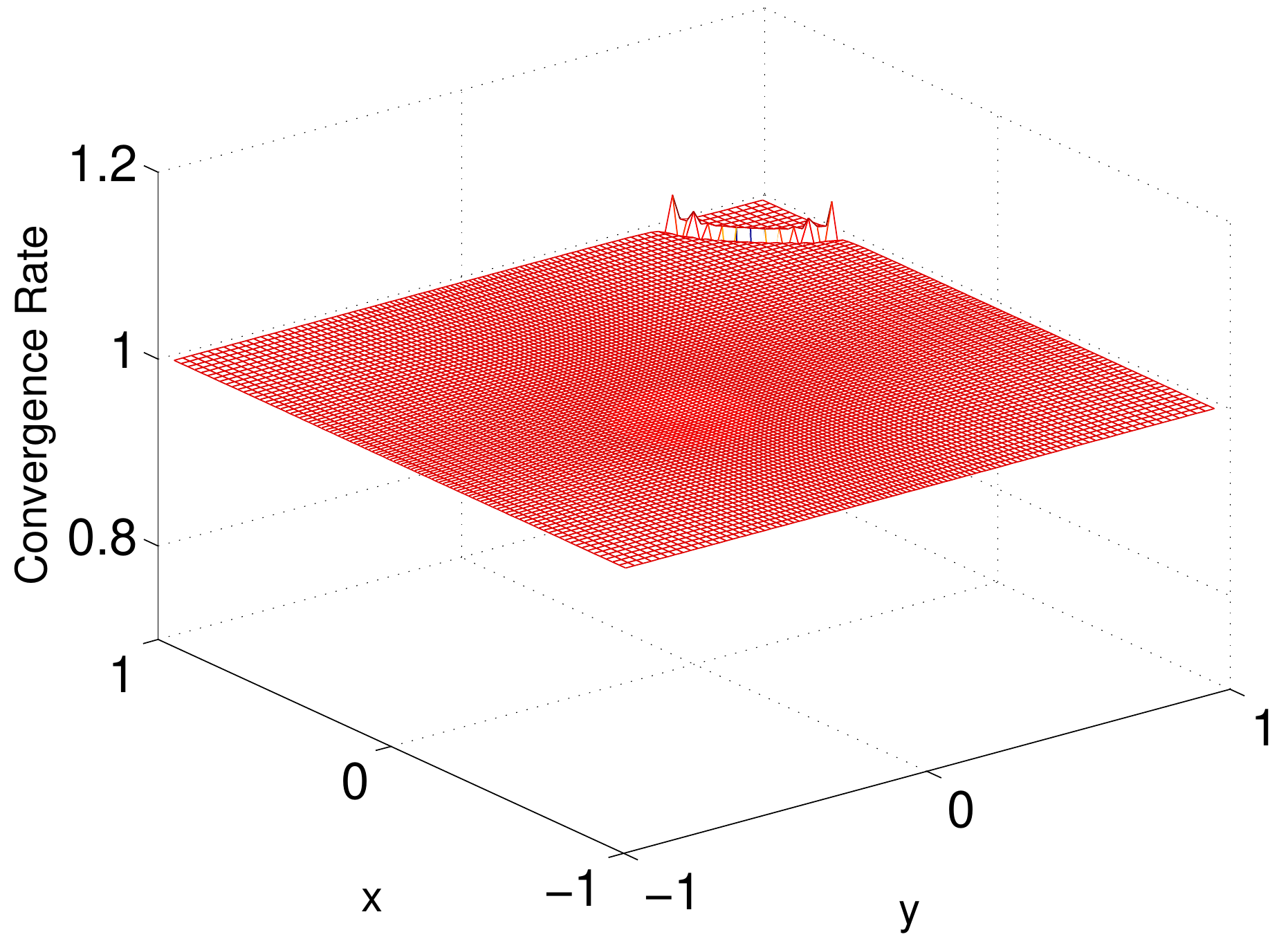,width=2.63in,height=1.68in}}  %%Code: Fig_7a_7b
{\epsfig{file=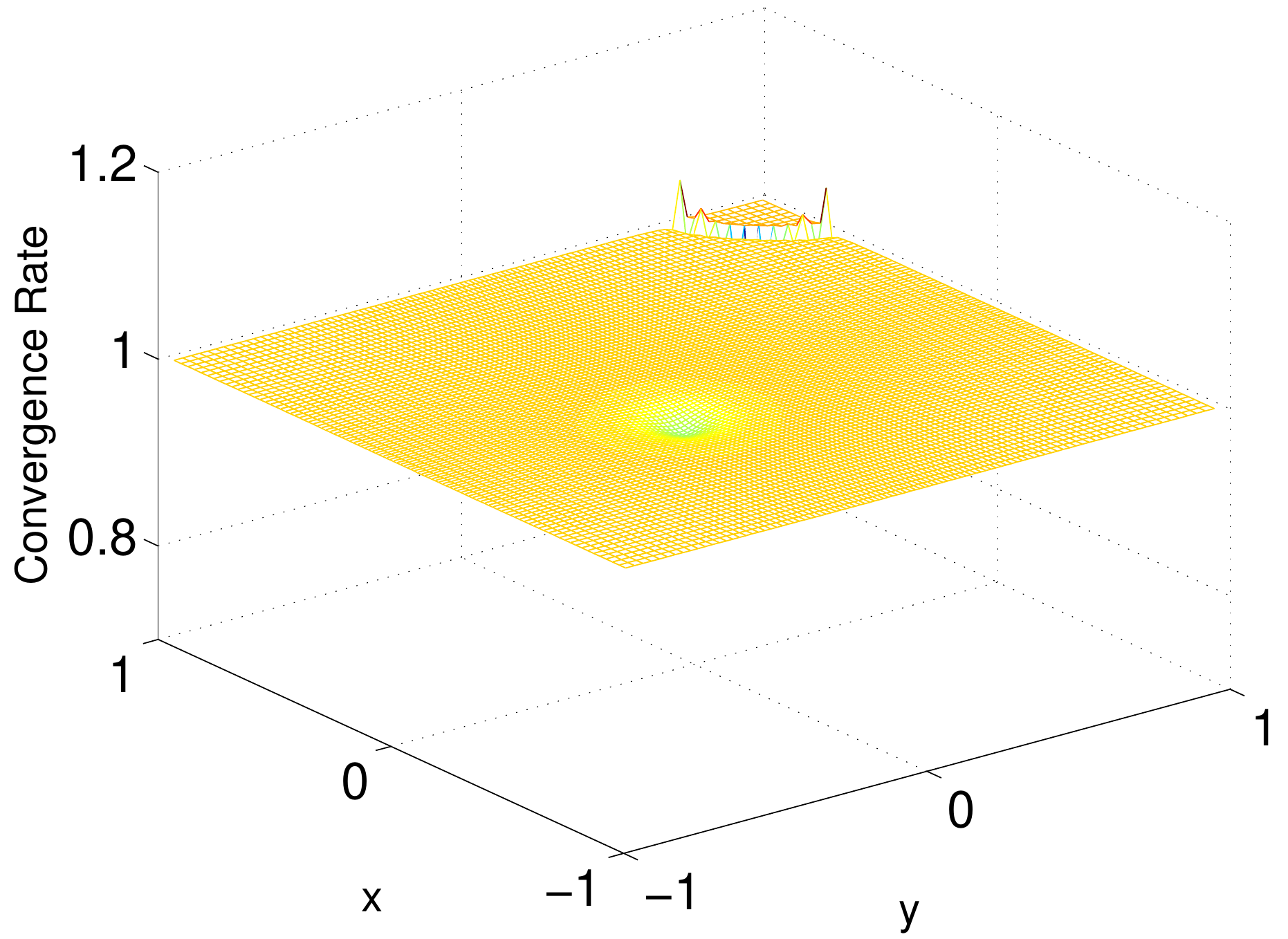,width=2.63in,height=1.68in}}  %%Code: Fig_7a_7b
\parbox[t]{12.8cm}{\scriptsize{\bf Figure 9.} Surface plots of the point-wise rates of convergence rates in time for the 
temporal derivative function of the solution, $u_t,$ at $T_1$ [LEFT] and $T_2$ [RIGHT] within the spatial domain. 
The mean convergence rates are approximately 0.991674641553278 and 0.986792663593675, respectively.}
\end{center}

As we have seen, surface illustrations will reflect the fact that the convergence issues can be restricted to a small 
number of mesh points near a corner of the domain due to the impacts from degeneracy. This is further recorded 
in Figure 9 which shows that the degeneracy has a noticeable effect on the temporal convergence of $u_t.$ There is an obvious ring of distinctive oscillatory waves near $(1,1)$ in the figure. However, these oscillations 
have little impact on the mean or median convergence rates which evolve from 0.991674641553278 to
0.986792663593675 in Figure 9. 
Further, the spectral convergence rates of $u_t$ evolve from 0.996338976463597 to 0.985326655425598.
These rates of convergence for the temporal derivative function are ideal, as they
remain approximately the same magnitudes as those from the previous example.

\section{Conclusions} \clearallnum

A semi-discretized method for solving degenerate stochastic Kawarada equations is derived and studied. 
The numerical method utilizes arbitrary grids in space and time. Various types of adaptations can
thus be introduced and incorporated to more precisely capture the underlying quenching, stochastic, and
degeneracy singularities. Exponential splitting strategies are employed for achieving a higher
effectiveness and efficiency in computations. The semi-discretized scheme is proven to preserve key physical
and mathematical properties such as the positivity and monotonicity of the solution under modest restrictions. 
The scheme acquired is stable in the classical von Neumann sense without additional constraints. An enriched 
stability analysis which incorporates the nonlinear stochastic quenching singularity is implemented. This extends 
the stability investigations in several recent publications \cite{Josh1,Josh2,Josh3}, and has proven to be more 
effective as compared with traditional analysis (see \cite{Beau1,Beau2,Beau3,Sheng5} and references therein). 

Since the Lax-Richtmeyer Equivalence theorem is not applicable to Kawarada equations which are highly
nonlinear, explorations of the nonlinear convergence of the aforementioned semi-discretized scheme are conducted through 
multiple numerical experiments. Our experiments indicate that, despite the strong quenching singularity, stochastic influence, and 
degeneracy, desirable rates of convergence can be expected. Although the degeneracy and quenching singularities do
affect the convergence, especially within certain regions, mean and median rates of convergence 
remain relatively unaffected. Our experiments also depict novel 
presentations of global convergence rates via convergence surfaces at various times. These surfaces 
provide important insights into the nature of nonlinear convergence. They also help to better understand 
how singular components in the differential equation may affect convergence as quenching is approached. 
%% Rates of convergence of the temporal derivative of the solution are also explored
%% ? in this paper. It is found the convergence rates behave similar to those of the solution. This is a surprising observation
%% ? since it is commonly understood that the derivative blows up as $t\rightarrow T_q.$ 

Our continuing endeavors include a study of the nonlinear stability and convergence of degenerate 
stochastic Kawarada equations. While the classical stability concept has been extended in our discussions, 
more rigorous nonlinear stability analysis and statistical tools will be implemented in future work. 
We intend to extend our study to Sobolev-like spaces and norms in the near
future, with the hopes that the current convergence approaches may have provided us with more valuable 
structural information. On the other hand, we also intend to %% extend the study of Kawarada quenching-combustion models by 
incorporate balanced nonlocal derivatives into the model \cite{Pag}. The inclusion of fractional partial differential 
equations may help to capture and utilize global features of the numerical combustion with solid applications to
the energy industry.




%\section*{References}

\bibliographystyle{elsarticle-num}
\bibliography{Dissertation_BibTex}


%\begin{thebibliography}{00}



%\end{thebibliography}
\end{document}